\documentclass[12pt]{article}
\usepackage {amsfonts,amsmath,amssymb}
\usepackage{hyperref}

\renewcommand{\a}{\mathfrak a}
\renewcommand{\c}{\mathfrak c}
\newcommand{\eps}{\varepsilon}
\newcommand{\g}{\mathfrak g}
\renewcommand{\l}{\mathfrak l}
\renewcommand{\t}{\mathfrak t}
\newcommand{\CC}{\mathcal C}
\newcommand{\C}{\mathbb C}
\newcommand{\CE}{\mathcal E}
\newcommand{\CH}{\mathcal H}
\newcommand{\CI}{\mathcal I}

\newcommand{\CL}{\mathcal L}
\newcommand{\R}{\mathbb R}
\newcommand{\CP}{\mathcal P}

\newcommand{\CS}{\mathcal S}
\newcommand{\CU}{\mathcal U}
\newcommand{\Z}{\mathbb Z}
\newcommand{\Sigmar}{\Sigma^\mathrm{r}}
\newcommand{\Hom}{\mathop{\rm Hom}\nolimits}

\newcommand{\tr}{\mathop{\rm tr}\nolimits}
\newcommand{\Ind}{\mathop{\rm Ind}\nolimits}
\renewcommand{\Re}{\mathop{\rm Re}\nolimits}
\renewcommand{\Im}{\mathop{\rm Im}\nolimits}
\newcommand{\vol}{\mathop{\rm vol}\nolimits}

\newcommand{\temp}{_{\mathrm{temp}}}
\newcommand{\disc}{_{\mathrm{disc}}}
\renewcommand{\ell}{_{\mathrm{ell}}}
\newcommand{\reg}{_{\mathrm{reg}}}

\newtheorem{theo}{Theorem}
\newtheorem{lemma}{Lemma}
\newtheorem{prop}{Proposition}
\newtheorem{cor}{Corollary}

\newcommand{\qed}{\hfill $\Box$\medskip}
\newenvironment{proof}[1][]{\par{\it Proof#1. }}{\qed\smallskip}

\begin{document}

\title{Asymptotic and descent formulas for weighted orbital integrals}
\author{Werner Hoffmann}
\maketitle
\begin{abstract}
We rewrite Arthur's asymptotic formula for weighted orbital integrals on real groups with the aid of a residue calculus and extend the resulting formula to the Schwartz space. Then we extract the available information about the coefficients in the decomposition of the Fourier transforms of Arthur's invariant distributions $I_M(\gamma)$ in terms of standard solutions of the pertinent holonomic system of differential equations. This allows us to determine some of those coefficients explicitly. Finally, we prove descent formulas for those differential equations, for their standard solutions and for the aforementioned coefficients, which reduce each of them to the case that $\gamma$ is elliptic in~$M$.

\end{abstract}

\section*{Introduction}

Weighted orbital integrals are the local contributions to the geometric side of the Arthur-Selberg trace formula. They are defined for connected reductive linear algebraic groups $G$ over a local field~$F$. A weighted orbital integral is a certain distribution supported on a conjugacy class $\{x^{-1}\gamma x\mid x\in G(F)\}$ in the group $G(F)$ of $F$-rational points and depends on a Levi $\R$-subgroup $M$ of $G$ containing~$\gamma$. Its value on a test function $f$ is denoted $J_M(\gamma,f)$, and it is tempered in the sense that it extends to continuous linear functionals on Harish-Chandra's Schwartz space~$\CC(G)$ of rapidly decreasing functions on~$G(\R)$. The precise definition will be recalled in the section~\ref{defs}. This paper is a contribution to the calculation of the Fourier transform of weighted orbital integrals. For a survey of this subject we refer to~\cite{ho-3}. It is important for certain applications of the trace formula, e.~g., for the determination of Gamma factors of zeta functions of Selberg type.

We first recall the notion of invariant Fourier transform. Let $\Pi\temp(G)$ the set of equivalence classes of irreducible tempered representations $\pi$ of~$G(\R)$. Temperedness means that
\[
\pi(f)=\int_{G(\R)}f(g)\pi(g)\,dg
\]
is well defined and bounded for $f\in\CC(G)$, where we have fixed a Haar measure $dg$ on~$G(\R)$. The character $\tr\pi(f)$ is then a tempered distribution (evaluated at the test function~$f$) that determines the class of~$\pi$. For fixed $f$, it can be viewed as a function $\hat f$ on~$\Pi\temp(G)$. We denote the image of the space $\CC(G)$ under this invariant Fourier transform by~$\CI(G)$. By definition, a tempered distribution $I:\CC(G)\to\C$ has the Fourier transform $\hat I:\CI(G)\to\C$ if
\[
\hat I(\hat f)=I(f)
\]
for all $f\in\CC(G)$. It is then necessarily invariant under inner automorphisms of~$G(\R)$.

The weighted orbital integral $J_M(\gamma,f)$ is non-invariant unless $M=G$, as are the weighted characters $\phi_M(\pi,f)$ for $\pi\in\Pi\temp(M)$ that are the local contributions to the spectral side of the trace formula. Starting from these two families of distributions, Arthur has constructed in~\cite{a-inv} tempered invariant distributions $I_M(\gamma,f)$ and used them to reformulate the trace formula in invariant form. It is these distributions for regular semisimple~$\gamma$ whose Fourier transforms we are actually interested in. Since our methods involve differential equations, we consider the groundfield~$F=\R$ only.

In order to describe the Fourier transform, one needs a description of the tempered dual $\Pi\temp(G)$ together with the space~$\CI(G)$ of functions on it. A Levi subgroup $L$ of~$G$ is by definition a Levi component of a parabolic subgroup $P$ of~$G$. Starting from a tempered representation $\sigma$ of $L(\R)$, one can produce a tempered representation of~$G(\R)$ using pull-back under the natural homomorphism $P\to L$ followed by induction. For generic~$\sigma$, this yields an irreducible representation $\sigma^G$. Every tempered representation is a direct summand of a representation induced from a square-integrable representation, and most of $\Pi\temp(G)$ is covered by the ranges of the sets $\Pi_2(L)$ of equivalence classes of square-integrable irreducible representations under induction. Moreover, the group of unramified tempered characters $L(\R)\to U(1)$ (i.~e., characters factoring through $\R_+$) acts freely on $\Pi_2(L)$ and provides it with the structure of a smooth manifold (in fact, a union of Euclidean spaces) such that $\hat f_L(\sigma)=\hat f(\sigma^G)$ extends to a Schwartz function invariant under a Weyl group~$W_L$.

The description of $\CI$ is complicated by the reducibility of induced representations. As a remedy, Arthur has introduced in~\cite{a-temp} the set $T\temp(G)$ consisting of virtual tempered characters, which are certain linear combinations of the (characters of the) constituents of a single induced representations. One defines induction $T\temp(L)\to T\temp(G)$ in the obvious way and defines the set $T\ell(G)$ of elliptic elements as those which are not induced from a proper Levi subgroup. This is again a union of Euclidean spaces, and $\CI$ can be identified with a direct sum of $W_L$-invariant Schwartz spaces on~$T\ell(L)$ over all Levi subgroups up to conjugacy.

It turns out that the Fourier transform of~$\hat I_M(\gamma)$ of $I_M(\gamma)$ is a measure whose atoms modulo~$i\a_G^*$ are always contained in a certain subset $T\disc(G)$ that can be larger than~$T\temp(G)$. Using the local trace formula, Arthur has proved in~\cite{a-Ftr} that
\[
\hat I_M(\gamma,\phi)=\sum_{[L]}
\int_{W_L\backslash T\disc(L)}\Phi_{M,L}(\gamma,\tau)\phi_L(\tau)\,d\tau
\]
for suitable smooth functions $\Phi_{M,L}$, where $\phi$ runs through~$\CI(G)$ and the measure $d\tau$ is invariant under~$i\a_L^*$.

Weighted orbital integrals as well as the distributions $I_M(\gamma)$ satisfy a system of linear differential equations as functions of $\gamma$ running through the regular points in a maximal torus $T$ of~$G$. This translates into a system of differential equations satisfied by the Fourier transforms $\Phi_{M,L}$. It has been proved in~\cite{ho-2} that there is a basis of standard solutions $\Psi_M^P$ on certain sectors depending on a parabolic subgroup~$P$ that are characterised by their asymptotic exponents $\mu\in\t^*$ as $\gamma\to\infty$. It follows that, roughly speaking,
\[
\Phi_{M,L}(\gamma,\tau)=\sum_{\mu}\sum_{S\supset M}c_{S,L}^P(\mu)\Psi_M^{P\cap S}(\gamma,\mu),
\]
where $S$ runs through the Levi subgroups containing~$M$ and $\mu$ runs through the parameters of the infinitesimal character of~$\tau$, which is assumed to be regular in a suitable sense. In order to determine the coefficients~$c$ as functions of $\mu$ and~$\tau$, one has to study the asymptotics of both sides as $\gamma$ tends to infinity and to the singular set in~$T(\R)$.

Arthur has proved in~\cite{a-asym} an asymptotic formula for the distributions $I_M(\gamma,f)$. This is somewhat subtle because just moving $\gamma$ would give zero tautologically. This one has to modify $f$ with the aid of a Schwartz multiplier at the same time. For Arthur's purposes of endoscopy, it was sufficient to work with test functions $f$ in the Hecke algebra. The answer was given by a contour integral in the complexified set $T(M)$ that avoids the poles of the intertwining operators on~$T\temp(M)$. He has sketched an argument to prove the existence of the limit for Schwartz functions $f\in\CC(G)$ that we take for granted.

Our aim in this paper is to convert the asymptotic formula into an expression that makes sense for Schwartz functions. This is done by shifting the contours of integration to the subset $T\temp(M)$ and taking care of the poles of the integrand by a suitable residue calculus. Next we study the asymptotics under the action of the multiplier on the expression in terms of $\Phi_{M,L}$. As a result, we get a formula for those coefficients $c_{S,L}^P(\mu)$ for which $\gamma^\mu$ is not decreasing on the sector corresponding to~$P$. In order to determine the remaining coefficients, one would have to know the asymptotics of the standard solutions $\Psi_M^P$ at the singular set, which are not yet available in general.

It is known that the distributions $I_M(\gamma)$ satisfy descent formulas, which reduce their study to the case that $\gamma$ is elliptic in~$M$. We prove similar descent formulas for the differential equations, their standard solutions~$\Psi_M^P$ and the coefficients~$c_{S,L}^P$. As an application, we obtain an explicit formula for $\Phi_{T,T}$, where $T$ is a split maximal torus, which was announced in~\cite{ho-3}.

This work was supported by the German Research Council (DFG), within the CRC 701.

\section{Weighted orbital integrals and invariant distributions}
\label{defs}
If $G$ is a connected reductive linear algebraic group defined over~$\R$, we write $\a_G$ for the real Lie algebra of the largest $\R$-split torus in the centre of~$G$. This is the only exception from our general rule to denote the complex Lie algebra of any algebraic group by the corresponding lower-case gothic letter. We have a natural homomorphism $H_G:G(\R)\to\a_G$ such that $H_G(\exp H)=H$ for $H\in\a_G$. For every $\lambda\in\a_{G,\C}^*$ and $\pi\in\Pi(G)$, we set $\pi_\lambda(g)=\pi(g)e^{\lambda(H_G(g))}$.

Every parabolic $\R$-subgroup $P$ of $G$ has a Levi decomposition $P=MN$, where $N$ is the unipotent radical of~$P$ and $M$ is a connected reductive $\R$-group just like~$G$. The set of roots of $\a_M$ in $N$ will be denoted by~$\Sigma_P$ and the subset of reduced roots by~$\Sigmar_P$. We fix a maximal compact subgroup $K$ of $G(\R)$; then $G(\R)=P(\R)K$. Writing $H_P(mnk)=H_M(m)$ for $m\in M(\R)$, $n\in N(\R)$ and $k\in K$, we obtain a continuous map $H_P:G(\R)\to\a_M$.

The set $\CP(M)$ of parabolic $\R$-subgroups $P$ with given Levi component~$M$ is in bijection with set of chambers $\a_P^+=\{H\in\a_M\,:\,\alpha(H)>0\mbox{ for $\alpha\in\Sigma_P$}\}$. For every $x\in G(\R)$, consider the following volume of a convex hull in $\a_M^G=\a_M/\a_G$:
\begin{equation}\label{vol}
v_M(x)=\mathop{\rm vol}\nolimits_{\a_M^G}\mathop{\rm conv}\{-H_P(x):P\in\CP(M)\}.
\end{equation}
This function is left $M$-invariant, because $H_P(mx)=H_M(m)+H_P(x)$.

  For $f$ in the Schwartz space $\CC(G)$ of rapidly decreasing $L^2$-functions on~$G(\R)$ and for $\gamma\in M(\R)$ such that the centralizer $G_\gamma$ of $\gamma$ in $G$ is contained in~$M$, the weighted orbital integral is defined as
\[
J_M(\gamma,f)=|D(\gamma)|^{1/2}\int_{G_\gamma(\R)\backslash G(\R)}f(x^{-1}\gamma x)v_M(x)\,d\dot x,
\]
where $D(\gamma)$ is the usual Weyl discriminant. Here we have fixed an invariant measure on $G_\gamma(\R)\backslash G(\R)$.

It has been shown~(\cite{a-inv}, Lemma~8.1) that the integral converges and defines a tempered distribution $J_M(\gamma)$ on~$G(\R)$, i.e., a continuous linear functional on~$\CC(G)$. Note that $v_G$ is constant equal to~1, so that $J_G(\gamma)$ is the ordinary (unweighted) orbital integral. This is the only case in which $J_M(\gamma)$ is invariant (under inner automorphisms).

The characteristic function of the convex hull in~(\ref{vol}) can be written as an alternating sum of characteristic functions of simplicial cones indexed by the elements of~$\CP(M)$. Their Fourier transforms as functions of $\lambda\in\a_{M,\C}^*$ converge for $\Re\lambda$ in a certain chamber of~$\a_M$, and in the limit one obtains
\begin{equation}\label{GM}
v_M(x)=\lim_{\lambda\to0}\sum_{P\in\CP(M)}\frac{e^{-\lambda(H_P(x))}}{\theta_P(\lambda)},
\end{equation}
where $\theta_P$ is the suitably normalized product of the linear functions defining the walls of~${}^+\a_P^*$, which is the dual cone of the chamber $\a_P^+$.

For $P\in\CP(M)$, $\pi\in\Pi\temp(M)$ and $\lambda\in\a_{M,\C}^*$, we realise the parabolically induced representation $\CI_P(\pi_\lambda)$ in the compact picture, i.~e. in a space of functions on our fixed maximal compact subgroup~$K$ that is independent of~$\lambda$. Given a second parabolic ${P'}\in\CP(M)$, the Knapp-Stein intertwining operators $J_{{P'}|P}(\pi_\lambda)$ from $\CI_P(\pi_\lambda)$ to~$\CI_{P'}(\pi_\lambda)$ is defined on $K$-finite vectors by a convergent integral for $\Re\lambda$ in a certain chamber and extends meromorphically to all~$\lambda$. The partial Plancherel densities $\mu_{{P'}|P}$ are defined by
\[
\mu_{{P'}|P}(\pi_\lambda)J_{P|{P'}}(\pi_\lambda)
J_{{P'}|P}(\pi_\lambda)=\mathrm{Id}.
\]
If $P$ and $P'$ are adjacent, i.~e.\ when $\Sigmar_{P'|P}=\Sigmar_{P'}\cap\Sigmar_{\bar P}$ consists of a single element~$\alpha$, we write $\mu_{P'|P}(\pi)=\mu_\alpha(\pi)$, so that in general
\begin{equation}\label{prodPlanch}
\mu_{P'|P}(\pi)=\prod_{\alpha\in\Sigmar_{P'|P}}\mu_\alpha(\pi).
\end{equation}
It is known that $\mu_\alpha(\pi_\lambda)$ as a function of~$\lambda\in\a_{M,\C}$ depends only on $\lambda(\check\alpha)$.

One can choose normalising factors $r_{{P'}|P}$ (see~\cite{a-Ftr}) such that the intertwining operators
\[
R_{{P'}|P}(\pi_\lambda)=r_{{P'}|P}(\pi_\lambda)^{-1}
J_{{P'}|P}(\pi_\lambda)
\]
have no poles for $\Re\lambda=0$ and satisfy, among other properties, an unrestricted transitivity. One has the $(G,M)$ family
\[
\mathcal R_{P'}(\lambda,\pi,P)=
R_{{P'}|P}(\pi)^{-1}R_{{P'}|P}(\pi_\lambda),\quad {P'}\in\CP(M),
\]
which in analogy to~(\ref{vol}) gives rise to the limit
\[
  \mathcal R_M(\pi,P)=\lim_{\lambda\to0}\sum_{{P'}\in\CP(M)}
\frac{\mathcal R_{P'}(\lambda,\pi,P)}{\theta_{P'}(\lambda)}.
\]
Although $\mathcal R_M(\pi,P)$ is a priori only defined on $K$-finite vectors, and the weighted character
\[
\phi_M^r(f,\pi)=\tr(\Ind_P^G(\pi,f)\mathcal R_M(\pi,P))
\]
only for $f$ in the Hecke algebra~$\CH(G)$, one can show (\cite{a-Ftr}, p.~175) that it extends to $f\in\CC(G)$ and provides a continuous map
\[
\phi_M^r:\CC(G)\to\CI(M).
\]
It follows from the transitivity and the intertwining property of $R_{{P'}|P}$ that $\phi_M^r$ is independent of~$P\in\CP(M)$, as the notation suggests. By definition, we have
\[
\phi_G^r(f)=\hat f.
\]

There is a canonical choice for the normalising factors $r_{{P'}|P}$ in terms of local $L$-factors (see~\cite{a-loc}), but in the analogous $p$-adic case they are not yet available. In order to make the normalisation universally canonic, rewrite
\[
\mathcal R_{P'}(\lambda,\pi,P)=
r_{P'}(\lambda,\pi,P)^{-1}\mathcal J_{P'}(\lambda,\pi,P)
\]
in terms of the $(G,M)$ families
\begin{align*}
\mathcal J_{P'}(\lambda,\pi,P)
&=J_{{P'}|P}(\pi)^{-1}J_{{P'}|P}(\pi_\lambda),\\
r_{P'}(\lambda,\pi,P)
&=r_{{P'}|P}(\pi)^{-1}r_{{P'}|P}(\pi_\lambda),
\end{align*}
which are defined as meromorphic functions on~$\a_{M,\C}^*$ for generic~$\pi$. Arthur has defined in~\cite{a-norm} the $(G,M)$ family
\[
m_{P'}(\lambda,\pi,P)=\mu_{{P'}|P}(\pi)\mu_{{P'}|P}(\pi_{\lambda/2})^{-1}
\]
that can replace $r_{P'}(\lambda,\pi,P)$ and gives rise to canonical weighted characters $\phi_M(f,\pi)$. Note that for these $(G,M)$ families we have limits $r_M(\pi,P)$ and $m_M(\pi,P)$ in analogy to~(\ref{vol}).

The invariant distributions that are the building blocks of the geometric side of the trace formula are characterised as follows (see~\cite{a-inv}, section~10, and~\cite{a-Ftr},~(3.5)). For all triples $(G,M,\gamma)$, where $G$ is a connected reductive $\R$-group, $M$ its Levi $\R$-subgroup and $\gamma\in M(\R)$ such that $G_\gamma\subset M$, there are invariant tempered distributions $I_M(\gamma)=I_M^G(\gamma)$ on $G(\R)$ such that
\[
J_M(\gamma,f)=\sum_{L\in\CL(M)}\hat I_M^L(\gamma,\phi_L(f)),
\]
where $\CL(M)$ denotes the set of all Levi subgroups of $G$ containing~$M$. Of course, there is also a version $I_M^r(\gamma)$ characterised by
\[
J_M(\gamma,f)=\sum_{L\in\CL(M)}\hat I_M^{r,L}(\gamma,\phi_L^r(f)).
\]

\section{The discrete part of the tempered dual}

Recall that the tempered dual $\Pi\temp(G)$ is exhausted by parabolically induced representations~$\CI_P(\sigma)$, where $P$ is a parabolic $\R$-subgroup with Levi component~$M$ and $\sigma\in\Pi_2(M)$. The description of the space~$\CI(G)$ of functions on $\Pi\temp(G)$ is complicated by the reducibility of~$\CI_P(\sigma)$ the for some~$\sigma$, that arises as follows.

If $w$ belongs to the stabiliser $W_\sigma$ of~$\sigma$ in~$W_M$, then the normalised intertwining operator $R_{wP|P}(\sigma)$ followed by translation with a representative of~$w$ yields a self-intertwining operator of~$\CI_P(\sigma)$. The elements $w$ for which that self-intertwining operator is a scalar multiple of the identity make up a normal subgroup $W_\sigma^0$ that is also related to the partial Plancherel densities in~(\ref{prodPlanch}). The set $\Sigmar_\sigma$ of those $\alpha\in\Sigmar(G,M)=\Sigmar_{\bar P|P}$ for which $\mu_\alpha(\sigma)=0$ is a root system with abstract Weyl group~$W_\sigma^0$. The reducibility of $\CI_P(\sigma)$ is measured by the $R$-group $R_\sigma=W_\sigma/W_\sigma^0$. One can fix a chamber $\c\subset\a_M$ for $\Sigmar_\sigma$ and identify $R_\sigma$ with the stabiliser of $\c$ in~$W_\sigma$.

In order to simplify the description of~$\CI(G)$, Arthur has introduced in~\cite{a-temp} a set $T(G)$ (and renamed it to $T\temp(G)$ in~\cite{a-asym}) as a replacement of $\Pi\temp(G)$. It consists of $G(\R)$-conjugacy classes $\tau$ of triples $(M,\sigma,r)$, where $r\in R_\sigma$. The discrete part $T\disc(G)$ mentioned in the heading consists of the classes of those triples for which the $W_\sigma^0$-coset of $r$ contains an element $w$ such that the fixed-point set $\a_M^w$ equals~$\a_G$. Often one does not distinguish between the class $\tau$ and the virtual character
\[
\hat f(\tau)=\tr(R(r)\CI_P(\sigma,f)),
\]
which is linear combination of (the characters of) the constituents of~$\CI_P(\sigma)$.

Parabolic induction from a Levi $\R$-subgroup $L$ extends to virtual representations and defines a map $T\temp(L)\to T\temp(G)$. In fact, if $\tau\in T\temp(L)$ is represented by the triple $(M,\sigma,r)$, where $r\in R_\sigma^L$, then the induced class $\tau^G\in T\temp(G)$ is represented by the same triple, where $r$ is now viewed as an element of the larger group $R_\sigma=R_\sigma^G$. The virtual character $\tau$ is elliptic if and only if $\a_M^r=\a_L$, and so $T\ell(L)\subset T\disc(L)$.

If $\tau_1\in T\ell(L_1)$ and $\tau_2\in T\ell(L_2)$ are such that $(L_1,\tau_1)$ and $(L_2,\tau_2)$ are $G(\R)$-conjugate, then clearly $\tau_1^G=\tau_2^G$. Actually, the converse is also true. Indeed, by Proposition~1.1 of~\cite{a-temp} we may assume that $\tau_i$ is represented by $(M,\sigma,r_i)$ with the same $M$ and~$\sigma$. Now the $r_i$ are $W_\sigma$-conjugate, and so are the $\a_{L_i}=\a_M^{r_i}$.

If $\tau\in T\temp(L)$ represented by a triple $(M,\sigma,r)$ and if $Q$, $Q'\in\CP(L)$, we set
\[
\mu_{Q'|Q}(\tau)=\mu_{Q'|Q}(\pi),\qquad r_{Q'|Q}(\tau)=r_{Q'|Q}(\pi),
\]
where $\pi$ is a constituent of~$\CI_R^L(\sigma)$, no matter which one. In fact,
\[
\mu_{Q'|Q}(\tau)=\mu_{RU'|RU}(\sigma),\qquad
r_{Q'|Q}(\tau)=r_{RU'|RU}(\sigma)
\]
for any $R\in\CP^L(M)$, where $U$ and $U'$ are the unipotent radicals of $Q$ and $Q'$.
These functions give rise to $(G,L)$ families and limits
\[
m_L(\tau,P)=m_L(\pi,P),\qquad r_L(\tau,P)=r_L(\pi,P).
\]
In analogy to~(\ref{prodPlanch}) we have the product formulae
\begin{equation}\label{prodPlancht}
\mu_{Q'|Q}(\tau)=\prod_{\beta\in\Sigmar_{Q'|Q}}\mu_\beta(\tau),\qquad
r_{Q'|Q}(\tau)=\prod_{\beta\in\Sigmar_{Q'|Q}}r_\beta(\tau),
\end{equation}
where $\mu_\beta(\tau)$ is the product of the $\mu_\alpha(\sigma)$ over those $\alpha\in\Sigmar_{RU'|RU}$ whose restriction to~$\a_L$ is a multiple of~$\beta$, and similarly in the case of the normalising factors.

For the rest of this section we assume that $\tau\in T\ell(L)$, again represented by a triple $(M,\sigma,r)$.
Let $\Sigmar_\tau$ be the set of those $\beta\in\Sigmar(G,L)$ for which  $\mu_\beta(\tau)=0$. The latter is the case if and only if there exists $\alpha\in\Sigmar_\sigma$ whose restriction to~$\a_L$ is a multiple of~$\beta$. Since $\mu_\alpha(\sigma)=\mu_{-\alpha}(\sigma)$, there is an even number $2n_\beta(\tau)$ of such roots~$\alpha$. As $\sigma$ is square-integrable, $\mu_\alpha(\sigma_\lambda)$ has a double zero at $\lambda(\check\alpha)=0$ for $\alpha\in\Sigmar_\sigma$. Thus, if $\beta\in\Sigmar_\tau$, the function $\mu_\beta(\tau_\lambda)$ has a zero of order $2n_\beta(\tau)$ and $r_\beta(\tau_\lambda)$ has a pole of order $n_\beta$ at $\lambda(\check\beta)=0$, where $\lambda$ runs through~$\a_{L,\C}^*$.

Later we shall need the constants
\begin{equation}\label{ntau}
k^L(\tau)=|R^L_{\sigma,r}|,\qquad
n(\tau)=n^G(\tau)=\sum_F \prod_{\beta\in F}\frac{n_\beta(\tau)}2,
\end{equation}
where $F$ runs through the bases of $(\a_L^G)^*$ contained in $\Sigmar_Q$ for some $Q\in\CP(L)$. Changing $Q$ entails replacing some $\beta$ by $-\beta$, but the value of $n_\beta(\tau)$ and hence $n(\tau)$ remains unchanged.
\begin{lemma}\label{Tdisc}
   The stabiliser $W_\tau$ of $\tau$ in $W_L$ acts transitively on the set of chambers of $\Sigmar_\tau$ in $\a_L$. In particular, it contains the reflections with respect to the elements of~$\Sigmar_\tau$. The virtual character $\tau^G$ belongs to $T\disc(G)$ if and only if $\Sigmar_\tau$ spans $(\a_L^G)^*$.
\end{lemma}

\begin{proof}
The chambers of $\Sigmar_\tau$ in $\a_L$ are of the form $\c\cap\a_L$, where $\c$ is an $r$-stable chamber of $\Sigmar_\sigma$ in~$\a_M$. Suppose that $\c'\subset\a_M$ is another chamber of this kind. Since $W_\sigma^0$ acts simply transitively on the set of chambers of $\Sigmar_\sigma$, it contains an element $w$ such that $\c=w\c'$. The element $wrw^{-1}r^{-1}$ of $W_\sigma$ belongs to the normal subgroup $W_\sigma^0$, and this element stabilises $\c$, so it must be the identity. It follows that $w$ commutes with $r$, so it stabilises $\a_L$ and restricts to an element of~$W_L$ which maps $\c\cap\a_L$ to $\c'\cap\a_L$. It is clear that $w$ fixes $\tau$, so the restriction belongs to~$W_\tau$. This shows the transitivity.

Given $\beta\in\Sigmar_\tau$, let $L_\beta\supset L$ be the Levi subgroup such that $\a_{L_\beta}$ is the kernel of~$\beta$, and let $P_\beta$, $P'_\beta$ be the two elements of~$\CP^{L_\beta}(L)$. Then $\mu_{P'_\beta|P_\beta}(\tau)=\mu_\beta(\tau)=0$ and hence $\beta\in\Sigma_\tau^{L_\beta,\mathrm r}$. By what we have proved, $W_\tau^{L_\beta}$ acts transitively on the set~$\{\a_{P_\beta}^+,\a_{P'_\beta}^+\}$, so it must contain the only nontrivial element of~$W_L^{L_\beta}$, viz.\ the reflection with respect to~$\beta$.

Let $H\in\a_L$ such that $\beta(H)=0$ for all $\beta\in\Sigmar_\tau$, i.~e. $\alpha(H)=0$ for all $\alpha\in\Sigmar_\sigma$. For all $w\in W_\sigma^0$ we have $rwH=rH=H$. If $\tau^G\in T\disc(G)$, then there exists $w\in W_\sigma^0$ such that $\a_M^{rw}=\a_G$, and it follows that $H\in\a_G$. This shows that $\Sigmar_\tau$ spans~$(\a_L^G)^*$.

Conversely, suppose the latter is the case. Then we can choose roots $\alpha_1$, \dots, $\alpha_n\in\Sigmar_\sigma$ whose restrictions to $\a_L^G$ form a basis of the dual space. Let $s_i$ denote the reflection with respect to $\alpha_i$ and set $w_0=s_1\cdots s_n\in W_\sigma^0$. It suffices to show that $\a_M^{rw_0}=\a_G$. Thus, let $H\in\a_M$ be a fixed point of~$rw_0$. Since the closure of $\c$ is a fundamental domain for the action of~$W_\sigma^0$, we can find $w$ in that group such that $H$ lies in $w\bar\c$. This chamber is stabilised by~$wrw^{-1}$, and after replacing $r$ and all the $s_i$ by their conjugates under $w$, $H$ by $wH$ and $\c$ by~$w\c$, we may assume that $rw_0$ fixes a point $H\in\bar\c$. Now $w_0H$ equals $r^{-1}H$, which lies in the fundamental domain $\bar\c$ of~$W_\sigma^0$, so that $H$ must be fixed by $w_0$. Applying $s_i$ to a point amounts to adding a multiple of the vector~$\check\alpha_i$, and applying $w_0$ amounts to adding a linear combination of these vectors, which is zero in case of the point~$H$. Since the $\check\alpha_i$ are linearly independent, each coefficient is zero, which means $s_i(H)=H$ and $\alpha_i(H)=0$ for all~$i$. But $H$ is also fixed by~$r$, hence in~$\a_L$, and by the choice of the $\alpha_i$ we have $H\in\a_G$, as was to be shown.
\end{proof}

\section{Shifting contours}

Let $M$ be a parabolic $\R$-subgroup of~$G$ and $P\in\CP(M)$. The asymptotic formula in \cite{a-asym} is stated as an integral over a subset
\[
T_\eps(M)=\{\tau_\eps\,:\,\tau\in T\temp(M)\}
\]
that involves the function $m_M(\tau,P)$ from the previous section. (Of course, the Levi subgroup in a triple representing $\tau \in T\temp(M)$ will in general be a proper subgroup of $M$ and has to be denoted by a different symbol.) We will shift the contour of integration to~$T\temp(M)$.

Recall that $T\temp(M)$ is the disjoint union of the sets $W_{L_1}^M\backslash T\ell(L_1)$ over all conjugacy classes of Levi subgroups $L_1$ of~$M$. More precisely, one should consider the images of the sets $T\ell(L_1)$ under the induction map $\tau\to\tau^M$. Let $T\ell(L_1)^1$ be the subset of those $\tau$ which are trivial on~$\exp\a_{L_1}$. Then every element of $T\ell(L_1)$ can be uniquely written as $\tau_\lambda$, where $\tau\in T\ell(L_1)^1$ and $\lambda\in i\a_{L_1}^*$.

\begin{lemma}\label{shift}
Let $\tau\in T\ell(L_1)^1$. Then for sufficiently small $\eps\in(\a_P^*)^+$ and every Paley-Wiener function $\phi$ on~$\a_{L_1,\C}^*$,
\begin{multline*}
\int_{i\a_{L_1}^*+\eps}\phi(\lambda)m_M(\tau_\lambda^M,P)\,d\lambda=\\
\sum_{L\in\CL(L_1)}\sum_{S\in\CL(M)}d_{L_1}^G(L,S)n^L(\tau)
\,\mathrm{p.v.}\!\int_{i\a_L^*}\phi(\lambda)m_M^S(\tau_\lambda^M,P\cap S)\,d\lambda.
\end{multline*}
The same formula is true with $m$ replaced by~$r$ on both sides.
\end{lemma}
The principal value means that we integrate over the complement of the $\delta$-neighbourhood of the union of root hyperplanes and let $\delta\to0$. Actually, we could restrict summation to pairs $(L,S)$ for which the constants $d_{L_1}^G(L,S)$ defined in \S7 of~\cite{a-invI} and $n^L(\tau)$ defined in~(\ref{ntau}) are nonzero. For such $L$ and all~$\lambda\in i\a_L^*$, we have $\tau_\lambda\in T\disc(L)$ by Lemma~\ref{Tdisc}.

\begin{proof}
The product formulae~(\ref{prodPlancht}) imply product formulae for the members of the corresponding $(G,L_1)$ families, viz.
\[
m_{Q'}(\lambda,\tau,Q)=\prod_{\beta\in\Sigmar_{Q'|Q}}
\frac{\mu_\beta(\tau)}{\mu_\beta(\tau_{\lambda/2})},\qquad
r_{Q'}(\lambda,\tau,Q)=\prod_{\beta\in\Sigmar_{Q'|Q}}
\frac{r_\beta(\tau_\lambda)}{r_\beta(\tau)}
\]
for $Q$, $Q'\in\CP(L_1)$ and $\lambda\in\a_{L_1,\C}^*$.
For every $\beta\in\Sigmar(A_{L_1},G)$, there are meromorphic functions $m_\beta'$ and $r_\beta'$ on $\C\times T\ell(L_1)^1$ such that
\[
-\frac{\partial({\textstyle\frac12}\beta)\mu_\beta(\tau_\lambda)}
{\mu_\beta(\tau_\lambda)}=m_\beta'(\lambda(\check\beta),\tau),\qquad
\frac{\partial(\beta)r_\beta(\tau_\lambda)}
{r_\beta(\tau_\lambda)}=r_\beta'(\lambda(\check\beta),\tau),
\]
where the directional derivative is defined as $\partial(\xi)\phi(\tau_\lambda)=\frac d{dt}\phi(\tau_{\lambda+t\xi})|_{t=0}$ for any $\xi\in \a_{L_1}^*$. Let us fix~$\tau$ and omit it from the notation. It follows from Harish-Chandra's explicit formula resp.\ the construction of normalising factors in terms of $L$-functions in~\cite{a-int}, (3.2) that $m_\beta'$ and $r_\beta'$ are functions of moderate growth uniformly in every vertical strip of bounded width except in a neighbourhood of a pole. Since $\tau$ is trivial on~$\exp\a_{L_1}$, the only possible pole on the imaginary axis is a simple one at zero with residue $-n_\beta(\tau)$. In the rest of the proof, we treat only the case of~$m$, that of $r$ being completely analogous.

If we fix a parabolic $Q_1\in\CP(L_1)$ contained in~$P$, then Corollary~7.3 in~\cite{a-EisII} (with $c_\beta(z)=1$ for $\beta\notin\Sigmar_{\bar Q_1}$) provides the formula
\begin{equation}\label{prodF}
m_M(\tau_\lambda^M,P)=\sum_F\vol(\a_M^G/\Z\check F_M)\prod_{\beta\in F}m_\beta'(\lambda(\check\beta)),
\end{equation}
where the sum is taken over all $F\subset\Sigmar_{\bar Q_1}$ such that $\check F_M=\{\check\beta_M\mid \beta\in F\}$ is a basis of $\a_M^G$. Here $\check\beta_M$ denotes the projection of the ``coroot'' $\check\beta\in\a_{L_1}$ to~$\a_M$.

For any linearly independent set $F\subset\Sigmar(A_{L_1},G)$, there is a unique Levi subgroup $L_F\in\CL(L_1)$ such that
\[
\a_{L_F}=\{H\in\a_{L_1}\mid\beta(H)=0\,\,\forall\beta\in F\}.
\]
If $\check F_M$ is a basis of $\a_M^G$, then $\a_M^G\oplus\a_{L_F}^G=\a_{L_1}^G$. We have fixed Haar measures on these vector spaces, and the constant $d_{L_1}^G(M,L_F)$ is such that for the corresponding Plancherel measures and any $\eps\in(\a_M^G)^*$ we have
\[
\int_{i(\a_{L_1}^G)^*+\eps}\phi(\lambda)\,d\lambda=d_{L_1}^G(M,L_F)\int_{i(\a_M^G)^*+\eps}\int_{i(\a_{L_F}^G)^*}\phi(\mu+\zeta)\,d\mu\,d\zeta.
\]
When we substitute $\phi(\lambda)m_M(\tau_\lambda^M,P)$ in place of~$\phi(\lambda)$, the second factor is independent of~$\mu$ and can be taken out of the inner integral. We use the isomorphism $(\a_M^G)_\C^*\to \C^F$ that maps $\zeta$ to the point $z$ with components $z_\beta=\zeta(\check\beta)$ and define a Paley-Wiener function on~$\C^F$ by
\[
\phi_F(z)=\int_{i(\a_{L_F}^G)^*}\phi(\mu+\zeta)\,d\mu.
\]
Since the elements of $F\subset\Sigmar_{\bar Q_1}$ are negative on~$\a_P^+$, the point corresponding to $\eps$ is of the form $-\eps_F$ with negative components~$-\eps_\beta$.

As $\zeta\in i(\a_M^G)^*$ is extended to a linear functional on~$\a_{L_1}^G$ vanishing on~$\a_{L_1}^M$, we have $\zeta(\check\beta)=\zeta(\check\beta_M)$. Thus our isomorphism is dual to the isomorphism $\R^F\to\a_M^G$ that maps the standard basis to $\check F_M$ and takes the Lebesgue measure to the Haar measure divided by $\vol(\a_M^G/\Z\check F_M)$. Recall that the Plancherel measure $\omega$ on~$i\R$ corresponding to the Haar measure on $\R$ is the Lebesgue measure divided by $2\pi$, and that Plancherel measures are inversely proportional to Haar measures. Thus, for generic $\eps\in(\a_P^*)^+$,
\[
\int_{i(\a_{L_1}^G)^*+\eps}\phi(\lambda)m_M(\tau_\lambda^M,P)\,d\lambda
=\sum_Fd_{L_1}^G(M,L_F)\int_{i\R^F-\eps_F}\phi_F(z)\prod_{\beta\in F}m_\beta'(z_\beta)\,\omega_F,
\]
where $\omega_F$ is the product measure of copies of $\omega$, or rather its translate.

For every Paley-Wiener function $\phi_\beta$ on $\C$ and $\eps_\beta>0$ sufficiently small we have
\[
\frac1{2\pi i}\int_{i\R-\eps_\beta}\phi_\beta(z)m_\beta'(z)\,dz=\frac{n_\beta(\tau)}2\,\phi_\beta(0)+\frac1{2\pi i}\,\mathrm{p.v.}\!\int_{i\R}\phi_\beta(z)m_\beta'(z)\,dz,
\]
where $i\R$ as a contour is endowed with the upward orientation. The differential form $dz/2\pi i$ induces the measure $\omega$ on the imaginary axis and its parallel translates. After coordinatewise application it follows that for every Paley-Wiener function $\phi_F$ on $\C^F$ and every $\eps_F\in\R^F$ with sufficiently small positive components,
\[
\int_{i\R^F-\eps_F}\phi_F(z)\prod_{\beta\in F}m_\beta'(z)\,\omega_F=\sum_{F',F''}\frac{n_{F'}(\tau)}{2^{|F'|}}\,\mathrm{p.v.}\!\int_{i\R^{F''}}\phi_F(z)\prod_{\beta\in F''}m_\beta'(z)\,\omega_{F''},
\]
where the sum is taken over all decompositions of $F$ into two disjoint subsets $F'$ and $F''$, and $n_{F'}=\prod_{\beta\in F'}n_\beta(\tau)$. On the right-hand side, we take principal values along the union of the coordinate hyperplanes in the obvious sense.

Next we have to take the sum over all sets $F$ as above. To every $F$ and its partition into $F'$ and $F''$, we can associate Levi subgroups $L_{F'}=L\in\CL(L_1)$ and $S\in\CL(M)$ defined by $\a_S=\a_{L_{F''}}\cap\a_M$. Then $\a_{L_1}^G=\a_L^G\oplus\a_S^G$, and the sets $F'\subset\Sigmar_{\bar Q_1}$ and $F_M''\subset\Sigmar_{\bar P}$ are bases of $(\a_{L_1}^L)^*$ and~$(\a_M^S)^*$, respectively. Conversely, any such configuration $L$, $S$, $F'$, $F''$ arises exactly once, hence we can replace the summation over $F$ by summations over $L$ and~$S$. The Hasse diagram below could be completed to a cube by the Levi subgroup corresponding to~$\a_M\cap\a_L$, but that does not turn up in our argument.
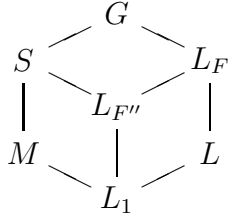
\begin{center}
\setlength{\unitlength}{1.5em}
\newcommand{\vl}{\makebox(1,1){\line(0,1){1}}}
\newcommand{\ul}{\makebox(1,1){\line(2,1){0.9}}}
\newcommand{\dl}{\makebox(1,1){\line(2,-1){0.9}}}
\begin{picture}(5,5)
\put(0,1){\makebox(1,1){$M$}}
\put(0,3){\makebox(1,1){$S$}}
\put(2,0){\makebox(1,1){$L_1$}}
\put(2,2){\makebox(1,1){$L_{F''}$}}
\put(2,4){\makebox(1,1){$G$}}
\put(4,1){\makebox(1,1){$L$}}
\put(4,3){\makebox(1,1){$L_F$}}
\put(0,2){\vl} \put(2,1){\vl} \put(4,2){\vl}
\put(1,0.5){\dl} \put(1,2.5){\dl} \put(3,3.5){\dl}
\put(1,3.5){\ul} \put(3,0.5){\ul} \put(3,2.5){\ul}
\end{picture}
\end{center}

When we re-substitute the defining expression of $\phi_F$ on the right-hand side, the argument of $\phi$ will lie in~$i(\a_{L}^G)^*$, and we have to make sense of $\omega_{F''}$ and the product of the $m_\beta'$ over~$\beta\in F''$. Like in the case of~$F$, we have an isomorphism $i\R^{F''}\to i(\a_M^S)^*$, under which the measure $\omega_{F''}$ corresponds to $\vol(\a_M^S/\Z\check F_M'')$ times the Plancherel measure. Considering the composition of natural isomorphisms $\a_{L_1}^M\to\a_{L_{F''}}^S\to\a_{L_F}^G$, we can write the coefficient in front of the integral as
\[
d_{L_1}^G(L_F,M)=d_{L_{F''}}^G(L_F,S)d_{L_1}^S(L_{F''},M).
\]
The second factor in this product times the Plancherel measure on~$i(\a_M^S)^*$ corresponds to the Plancherel measure on $i(\a_{L_1}^{L_{F''}})^*$, which in turn corresponds to $d_{L_1}^{L_F}(L,L_{F''})$ times the Plancherel measure on~$i(\a_L^{L_F})^*$. We can now combine both integrals into one integral over~$i(\a_L^G)^*$. The resulting coefficient
\[
d_{L_1}^{L_F}(L,L_{F''})d_{L_{F''}}^G(L_F,S)=d_{L_1}^G(L,S)
\]
can be taken out of the sum over~$F''$, which becomes
\[
\sum_{F''}\vol(\a_M^S/\Z\check F_M'')\prod_{\beta\in F''}m_\beta'(\lambda(\check\beta))=m_M^S(\tau_\lambda^M,P\cap S).
\]
The sum over $F'$ gives rise to the coefficient~$n^L(\tau)$. Finally, we replace $\phi$ by a Paley-Wiener function on $i\a_{L_1}^*$ averaged over~$i\a_G^*$.
\end{proof}

\section{Temperedness}

In \S7 of~\cite{a-asym}, Arthur has sketched a proof that the limit in the asymptotic formula exists uniformly for test functions $f$ in the Schwartz space $\CC(G)$ and that the resulting distribution is tempered. However, that distribution is given by an integral over the set $T_\eps(M)$ where the Fourier transforms $\hat f_M$ of such functions are not defined. We are going to show that after shifting contours to the set $T\temp(M)$, the resulting expression does converge for Schwartz functions.

This requires certain bounds in terms of infinitesimal characters of the elements of~$T\temp(M)$, which we write again as $\tau^M$, where $\tau\in T\ell(L_1)$. Suppose that $\tau$ is the $L_1$-conjugacy class of $(M_1,\sigma,r)$ and that $T_1$ is a maximal torus in~$M_1$. Then the infinitesimal character of~$\sigma$ is parametrised by a $W_{\t_1}^{\mathfrak m_1}$-orbit in $\t_1^*$ that we denote~$\t_1^*(\sigma)$. The infinitesimal character of the virtual representation $\tau$ is that of any of its components. It is parametrised the $W_{\mathfrak t_1}^{\l_1}$-orbit $\t_1^*(\tau)$ containing~$\t_1^*(\sigma)$, and similarly $\t_1^*(\tau^M)$ is the $W_{\mathfrak t_1}^{\mathfrak m}$-orbit containing~$\t_1^*(\sigma)$. We assume that a $W_{\mathfrak t_1}^\g$-invariant inner product has been fixed on~$\mathfrak \t_1^*$ compatible with all the decompositions $\mathfrak \t_1=\mathfrak \t_1^{L_1}\oplus\a_{L_1,\C}$. If $\tau\in T\ell(L_1)^1$ and $\lambda\in i\a_{L_1}^*$, then
\[
\t_1^*(\tau_\lambda)=\{\mu+\lambda\mid \mu\in\t_1^*(\tau)\}.
\]

We say that a function $\phi$ on $T\ell(L_1)$ is smooth and slowly increasing if for every differential operator $D$ with constant coefficients on $i\a_{L_1}^*$ there exist a natural number $n$ and a constant $c>0$ such that
\[
|D\phi(\tau)|\le c(1+\|\mu\|)^n\qquad\mbox{for $\mu\in\t^*(\tau)$.}
\]

\begin{lemma}\label{tempext}
For all Levi subgroups $L_1\subset M$ of $G$, parabolics $P\in\CP(M)$ and smooth slowly increasing functions $\phi$ on $T\ell(L_1)$, the function
\[
\sum_{w\in W_\tau}
\phi(\tau_{w\lambda})m_M(\tau_{w\lambda}^M,P)
\]
extends to a smooth slowly increasing function on $T\ell(L_1)^1\times i\a_{L_1}^*$. The same is true with $m$ replaced by~$r$.
\end{lemma}
Recall that $m_M(\tau_\lambda^M,P)$ and $r_M(\tau_\lambda^M,P)$ are meromorphic functions of $\lambda\in\a_{L_1,\C}^*$ with poles along the hyperplanes $\lambda(\check\beta)=0$ for each $\beta\in\Sigmar_\tau$.
\begin{proof}
Due to formula~(\ref{prodF}) and its analogue, it suffices to consider the functions
\[
\prod_{\beta\in F}m_\beta'(\lambda(\check\beta),\tau),\qquad
\prod_{\beta\in F}r_\beta'(\lambda(\check\beta),\tau)
\]
in place of~$m_M(\tau_\lambda^M,P)$ and $r_M(\tau_\lambda^M,P)$, where $F\subset\Sigmar_\tau$ is linearly independent. It follows from Harish-Chandra's explicit formula resp.\ the construction of normalising factors in terms of $L$-functions in~\cite{a-int}, (3.2) that
\[
m_\beta'(\lambda(\check\beta),\tau)+\frac{n_\beta(\tau)}{\lambda(\check\beta)},\qquad
r_\beta'(\lambda(\check\beta),\tau)+\frac{n_\beta(\tau)}{\lambda(\check\beta)}
\]
extend to smooth slowly increasing functions of~$\tau_\lambda$.
Thus we need only consider the reciprocal of
\[
\Pi_F(\lambda)=\prod_{\beta\in F}\lambda(\check\beta)
\]
in place of~$m_M(\tau_\lambda^M,P)$ and $r_M(\tau_\lambda^M,P)$.

We choose a chamber $\mathfrak c$ of $\Sigmar_\tau$, denote by $\Sigmar_{\mathfrak c}$ the set of those elements of $\Sigmar_\tau$ which are positive on~$\mathfrak c$ and set
\[
\Pi_{\mathfrak c}(\lambda)=\prod_{\beta\in\Sigmar_{\mathfrak c}}\lambda(\check\beta).
\]
Note that $\Pi_{\mathfrak c}(w\lambda)=\eps_\tau(w)\Pi_{\mathfrak c}(\lambda)$, where $\eps_\tau:W_\tau\to\{\pm1\}$ is a character independent of the choice of~$\mathfrak c$. If we choose $\mathfrak c$ so that the kernel of a given $\beta\in\Sigmar_\tau$ is one of its walls, then we get for the reflection $s$ corresponding to~$\beta$ that $\eps_\tau(s)=-1$.

Our expression becomes
\[
\sum_{w\in W_\tau}\frac{\phi(\tau_{w\lambda})}{\Pi_F(w\lambda)}
=\frac1{\Pi_{\mathfrak c}(\lambda)}
\sum_{w\in W_\tau}\eps_\tau(w)\phi(\tau_{w\lambda})
\frac{\Pi_{\mathfrak c}(w\lambda)}{\Pi_F(w\lambda)}.
\]
Since $\Pi_{\mathfrak c}/\Pi_F$ extends to a polynomial, every term in the sum extends to a smooth slowly increasing function. If $s$ is the reflection corresponding to~$\beta\in\Sigmar_\tau$, then the contributions of $w$ and $ws$ to the sum on the right-hand side cancel for $s\lambda=\lambda$. Thus the whole sum vanishes for $\lambda(\check\beta)=0$.

It follows from Taylor's theorem that for every smooth slowly increasing function $\phi_0$ that vanishes for $\lambda(\check\beta)=0$ there is a smooth slowly increasing function $\phi_1$ such that
\[
\phi_0(\tau_\lambda)=\lambda(\check\beta)\phi_1(\tau_\lambda)
\]
and that
\[
|\phi_1(\tau_\lambda)|
\le\sup_{|t|\le1}\left(|\phi_0(\tau_\lambda)|
+|\partial(\check\beta)\phi_0(\tau_{\lambda+t\check\beta})|\right).
\]
Applying this repeatedly to the last sum over $W_\tau$, we see that it is divisible by $\Pi_{\mathfrak c}(\lambda)$ as a smooth slowly increasing function.
\end{proof}

\section{The asymptotic formula}

Roughly speaking, the asymptotic formula of~\cite{a-asym} describes the limit of $I_M(\gamma_X)$, where $\gamma\in M\reg(\R)$ and $\gamma_X=\gamma\exp X$, as $X$ tends to infinity in a certain chamber of~$\a_M$. More precisely, one fixes $P\in\CP(M)$ and $r>0$ and lets $\|X\|$ tend to infinity while $\alpha(X)>r\|X\|$ for every fundamental root $\alpha$ of $\a_M$ in the unipotent radical of~$P$. This is abbreviated as $X\mathrel{\mathop{\longrightarrow}\limits_{P,r}}\infty$. In order to get a limit that is not trivially zero, one has to replace the test function $f\in\CC(G)$ by a new function $f_X\in\CC(G)$ that varies with $X$ and is characterised with the help of a multiplier $\alpha_X$ as follows.

If $M_1$ is a Levi subgroup of~$G$, $\sigma\in\Pi_2(M_1)$ and $P_1\in\CP(M_1)$, then $\CI_{P_1}(\sigma,f_X)=0$ unless a conjugate of $M_1$ is contained in~$M$, and if $M_1\subset M$, then
\[
\CI_{P_1}(\sigma,f_X)=\hat\alpha_X(\sigma)\CI_{P_1}(\sigma,f),
\]
where
\[
\hat\alpha_X(\sigma)=\frac1{|W_{M_1}|}\sum_{w\in W_{M_1}}e^{w\nu(X)}.
\]
Here $\nu\in i\a_{M_1}^*$ denotes the infinitesimal central character of~$\sigma$ defined by $\sigma(\exp H)=e^{\nu(H)}\mathrm{Id}$ for $H\in\a_{M_1}$. Note that for $M=G$ one simply gets $f_X(\gamma_X)=f(\gamma)$.

If $T_1$ is a maximal torus in~$M_1$, then the lattice of cocharacters of~$T_1$ over~$\C$ can be naturally embedded into~$\t_1$, and for any element $Z$ of its $\R$-span and any $\mu\in\t_1^*(\sigma)$ we have $\nu(Z)=\Im\mu(Z)$. That is why in~\cite{a-asym} the infinitesimal central character was called the imaginary part of the infinitesimal character.

For a virtual character $\tau\in T\temp(G)$ represented by a triple $(M_1,\sigma,r)$, where $\sigma\in\Pi_2(M_1)$ and $r\in R_\sigma$, one writes
\[
\hat f(\tau)=\tr(R(r,\sigma)\CI_{P_1}(\sigma,f)).
\]
It is clear that
\[
(f_X)_G(\tau)=\hat\alpha_X(\tau)\hat f(\tau),
\]
where $\hat\alpha_X(\tau)=\hat\alpha_X(\sigma)$ if a conjugate of $M_1$ is contained in~$M$ and $\hat\alpha_X(\tau)=0$ otherwise.

More generally, let $L$ be a Levi subgroup containing~$M_1$, and let $\tau\in T\temp(L)$ be the $L(\R)$-conjugacy class of the triple $(M_1,\sigma,r)$, where $r\in R_\sigma^L$. Then one uses the notation
\[
\hat f_L(\tau)=\hat f(\tau^G),
\]
and by transitivity of induction we get
\[
\widehat{(f_X)}_L(\tau)=\hat\alpha_X(\tau)\hat f_L(\tau),
\]
with $\hat\alpha_X(\tau)=\hat\alpha_X(\tau^G)$ as above.

As with all objects that depend on Levi subgroups up to conjugacy, it is sufficient to consider only Levi subgroups in $\CL=\CL(M_0)$ for a fixed minimal Levi $\R$-subgroup~$M_0$. We write $W_0^G$ for $W_{M_0}^G$.

The general structure of the Fourier transforms of the invariant distributions $I_M(\gamma)$ is given by Theorem~4.1 of~\cite{a-Ftr}. It asserts that there are unique smooth functions $\Phi_{M,L}$ on $(M(\R)\cap G\reg(\R))\times T\disc(L)$ such that, for all $f\in\CC(G)$,
\begin{equation}\label{Ftr}
I_M(\gamma,f)=\sum_{L\in\CL}\frac{|W_0^L|}{|W_0^G|}
\int_{T\disc(L)}\Phi_{M,L}(\gamma,\check\tau)\hat f_L(\tau)\,d\tau
\end{equation}
and that $\Phi_{M,L}(\gamma,\tau)$ is invariant under $W_0^G$ acting on the pairs $(L,\tau)$. Moreover, those functions are slowly increasing in the variable~$\tau$, which ensures the absolute convergence of the integrals. In the special case that $M=G$, the functions $\Phi_{G,L}(\gamma,\tau)$ vanish unless $\tau\in T\ell(L)$. We have written the contragredient $\check\tau$ of~$\tau$ as an argument of~$\Phi_{M,L}$ in order to simplify formulas in the next section.
\begin{theo}\label{asymI}
  Let $P\in\CP(M)$, $\gamma\in M\reg(\R)$ and $f\in\CC(G)$. Then
\begin{multline*}
\lim\limits_{X\mathrel{\mathop{\longrightarrow}\limits_{P,r}}\infty}
I_M(\gamma_X,f_X)=
\sum_{\substack{L_1,L\in\CL\\L_1\subset L}}
\frac{|W_0^{L_1}|}{|W_0^G|}
\sum_{\tau\in T\ell(L_1)^1}
\int_{i\a_L^*}\hat f_{L_1}(\tau_\lambda)\frac{n^L(\tau)}{k^{L_1}(\tau)}\\
\times \sum_{S\in\CL(L_1)}d_{L_1}^G(L,S)
\sum_{\substack{w\in W_0^G/W_0^M\\L_1\subset wM\subset S}}
\Phi_{wM,L_1}^{wM}(w\gamma,\check\tau_{-\lambda})
m_{wM}^S(\tau_\lambda^{wM},wP\cap S)\,d\lambda,
\end{multline*}
where $n^L(\tau)$ and $k^{L_1}(\tau)$ are defined in~(\ref{ntau}). The integrand is $W_0^G$-invariant as a function of $(L_1,L,\tau,\lambda)$, and the integral-sum over $\tau$ and $\lambda$ is absolutely convergent.

For every family $r_{P'|P}$ of normalising factors as in~\cite{a-Ftr}, section~2, the analogous formula for the limit of $I_M^r(\gamma_X,f_X)$ with $m_M^S$ replaced by~$r_M^S$ is valid.
\end{theo}
\begin{proof}
  The special case of equation~(\ref{Ftr}) with $G$ replaced by~$M$ reads
\[
I_M^M(\gamma,h)=\sum_{L_1\in\CL^M}\frac{|W_0^{L_1}|}{|W_0^M|}
\int_{T\ell(L_1)}\Phi_{M,L_1}^M(\gamma,\check\tau)h_{L_1}(\tau)\,d\tau
\]
for $h\in\CC(M)$. For $\lambda\in i\a_M^*$ we have
\[
\Phi_{M,L_1}^M(\gamma,\tau_\lambda)
=\Phi_{M,L_1}^M(\gamma,\tau)e^{\lambda(H_M(\gamma))}.
\]
Following \cite{a-asym}, \S5, we use this equation in the case $\lambda\in\a_{M,\C}^*$ to extend the definition of~$\Phi_{M,L_1}^M$. If $h\in\CH(M)$, so that $h_{L_1}$ is a Paley-Wiener function, this allows one to replace $\tau$ by $\tau_\eps$, which results in a shift of the contour of integration like in Lemma~\ref{shift}. As in that lemma, we write the elements of $T\ell(L_1)$ in the form $\tau_\lambda$ with $\tau\in T\ell(L_1)^1$ and $\lambda\in i\a_{L_1}$. According to \cite{a-temp}, equation~(3.5), the restriction of the measure on $T\ell(L_1)$ to an $i\a_{L_1}$-orbit is the image of the measure on $i\a_{L_1}$ divided by~$k^{L_1}(\tau)$.

For the time being, let $f\in\CH(G)$. Note that $\gamma_X\in G\reg(\R)$ for $X\in\a_P^+$ sufficiently far away from the walls. According to Corollary~6.2 of~\cite{a-asym}, for $\eps\in(\a_P^*)^+$ small enough, the left-hand side of the equation in the theorem equals
\begin{equation}\label{rhslim}
\sum_{L_1\in\CL^M}\frac{|W_0^{L_1}|}{|W_0^M|}
\sum_{\tau\in T\ell(L_1)^1}
k^{L_1}(\tau)^{-1}\int_{i\a_{L_1}^*\!+\eps}
\Phi_{M,L_1}^M(\gamma,\check\tau_{-\lambda})
\hat f_{L_1}(\tau_\lambda)m_M(\tau_\lambda^M,P)\,d\lambda.
\end{equation}
If we apply Lemma~\ref{shift}, this expression becomes
\begin{multline*}
\sum_{L_1\in\CL^M}\frac{|W_0^{L_1}|}{|W_0^M|}
\sum_{\substack{L\in\CL(L_1)\\S\in\CL(M)}}
\sum_{\tau\in T\ell(L_1)^1}
d_{L_1}^G(L,S)\frac{n^L(\tau)}{k^{L_1}(\tau)}\\
\mathrm{p.v.}\!\int_{i\a_L^*}
\hat f_{L_1}(\tau_\lambda)
\Phi_{M,L_1}^M(\gamma,\check\tau_{-\lambda})
m_M^S(\tau_\lambda^M,P\cap S)\,d\lambda.
\end{multline*}
The value of this expression remains unchanged if we replace all occurrences of $M$, $P$, $\gamma$, $L_1$, $L$, $S$, $\tau$ and $\lambda$ by their conjugates under a fixed~$w\in W_0^G$. The same is true if we sum the resulting expression over $w$ and divide by~$|W_0^G|$. Now we denote the conjugated objects $wL_1$, $wL$, $wS$, $w\tau$ and $w\lambda$ by the original variables and obtain
\begin{multline*}
\frac1{|W_0^G|}\sum_{w\in W_0^G}
\sum_{L_1\in\CL^{wM}}\frac{|W_0^{L_1}|}{|W_0^M|}
\sum_{\substack{L\in\CL(L_1)\\S\in\CL(wM)}}
\sum_{\tau\in T\ell(L_1)^1}
d_{L_1}^G(L,S)\frac{n^L(\tau)}{k^{L_1}(\tau)}\\
\mathrm{p.v.}\!\int_{i\a_L^*}
\hat f_{L_1}(\tau_\lambda)
\Phi_{wM,L_1}^{wM}(w\gamma,\check\tau_{-\lambda})
m_{wM}^S(\tau_\lambda^{wM},{wP}\cap S)\,d\lambda.
\end{multline*}
We change the order of summation and integration to get
\begin{multline*}
\sum_{\substack{L_1,L\in\CL\\L_1\subset L}}
\frac{|W_0^{L_1}|}{|W_0^G|}
\sum_{\tau\in T\ell(L_1)^1}
\mathrm{p.v.}\!\int_{i\a_L^*}\hat f_{L_1}(\tau_\lambda)
\frac{n^L(\tau)}{k^{L_1}(\tau)}\sum_{S\in\CL(L_1)}d_{L_1}^G(L,S)\\
\times\frac1{|W_0^M|}\sum_{\substack{w\in W_0^G\\L_1\subset wM\subset S}}
\Phi_{wM,L_1}^{wM}(w\gamma,\check\tau_{-\lambda})
m_{wM}^S(\tau_\lambda^{wM},wP\cap S)\,d\lambda.
\end{multline*}
Note that $\Phi_{M,L_1}^M(\gamma,\tau)$ depends only on the $M(\R)$-conjugacy class of $\gamma$, so the term corresponding to $w$ depends only on the right coset of $w$ modulo~$W_0^M$.

In \S7 of~\cite{a-asym}, Arthur has sketched a proof that the limit $X\mathrel{\mathop{\longrightarrow}\limits_{P,r}}\infty$ also exists for $f\in\CC(G)$ and is a tempered distribution. Assuming that this has been carried out, it remains to show that the above integral-sum converges not only as a principal value, but absolutely and defines a tempered distribution.

Substituting $w^{-1}$ for~$w$, we can write the sum over $W_0^G$ as
\[
\sum_{\substack{w\in W_0^G\\wL_1\subset M\subset wS}}
\Phi_{M,wL_1}^{M}(\gamma,(w\check\tau)_{-w\lambda})
m_M^{wS}((w\tau)_{w\lambda}^M,P\cap wS).
\]
We consider the partial sum over those $w$ for which $wL_1$ equals a fixed Levi subgroup~$L_1'\subset M$. These make up a left coset modulo the stabiliser of~$L_1'$ in $W_0^G$, whose quotient modulo~$W_0^{L_1'}$ equals~$W_{L_1'}$. We further restrict summation to the elements $w\in W_{L_1'}$ for which $wS$ equals a fixed Levi subgroup~$S'\supset M$, and among those to a left coset modulo~$W_{L_1'}^{S'}$. Finally, we restrict the sum to those $w$ for which $w\tau$ equals a fixed~$\tau'\in T\ell(L_1')^1$. They make up a left coset modulo~$W_{\tau'}^{S'}$, and our partial sum becomes
\[
|W_0^{L_1'}|\sum_{w'\in W_{\tau'}^{S'}}
\Phi_{M,L_1'}^M(\gamma,\check\tau'_{-w'\lambda'})
m_M^{S'}(\tau_{w'\lambda'}'^M,P\cap S'),
\]
where $\lambda'=w\lambda$. Now the convergence of the integral-sum and the continuous dependence on $f\in\CC(G)$ follow from Lemma~\ref{tempext} applied to $S'$ in place of~$G$.

Note that Corollary~6.2 of~\cite{a-asym} has an obvious analogue for $I_M^r(\gamma)$ with $m_M(\tau,P)$ replaced by~$r_M(\tau,P)$ that follows from the same Theorem~6.1. The corresponding version of our Theorem can thus be proved in a parallel fashion.
\end{proof}

\begin{cor}\label{PhiP}
Let $P\in\CP(M)$. Then there are unique smooth functions $\Phi_{P,L}$ on $M\reg(\R)\times T\disc(L)$ such that, for all $f\in\CC(G)$,
\[
\lim_{X\mathrel{\mathop{\longrightarrow}\limits_{P,r}}\infty}I_M(\gamma_X,f_X)
=\sum_{L\in\CL}\frac{|W_0^L|}{|W_0^G|}
\int_{T\disc(L)}\Phi_{P,L}(\gamma,\check\tau)\hat f_L(\tau)\,d\tau
\]
and that $\Phi_{P,L}(\gamma,\tau)$ is invariant under $W_0^G$ acting on the pairs $(L,\tau)$. Every element of $T\disc(L)$ is of the form $\tau^L$ for some Levi subgroup $L_1$ of $L$ and some $\tau\in T\ell(L_1)$, and we have
\begin{multline*}
\Phi_{P,L}(\gamma,\tau^L)
=k_{L_1}^L(\tau)n^L(\tau)\sum_{S\in\CL(L_1)}d_{L_1}^G(L,S)\\
\times\sum_{\substack{w\in W_0^G/W_0^M\\L_1\subset wM\subset S}}
\Phi_{wM,L_1}^{wM}(w\gamma,\tau)
m_{wM}^S(\check\tau^{wM},wP\cap S),
\end{multline*}
where $k_{L_1}^L(\tau)=k^L(\tau)/k^{L_1}(\tau)$. In particular,  $\Phi_{P,L}(\gamma,\tau^L)=0$ unless a conjugate of $L_1$ is contained in~$M$.
\end{cor}

The formula is an easy consequence of Theorem~\ref{asymI} and the definition
\[
\int_{T\disc(L)}\phi(\tau)\,d\tau
=\sum_{L_1\in\CL^L}\frac{|W_0^{L_1}|}{|W_0^L|}
\sum_{\substack{\tau\in T\ell(L_1)^1\\n^L(\tau)\ne0}}k^L(\tau)^{-1}\int_{i\a_L^*}\phi(\tau_\lambda^L)\,d\lambda,
\]
of the measure on $T\disc(G)$, which is equivalent to equation~(3.5) of~\cite{a-temp}. The smoothness of $\Phi_{P,L}$ follows from Lemma~\ref{tempext} by the argument in the proof of Theorem~\ref{asymI}, and this implies the uniqueness in view of the trace Paley-Wiener theorem in~\cite{a-PW}.

Obviously, Corollary~\ref{PhiP} has an analogue for the distribution $I_M^r$, its Fourier transform $\Phi_{M,L}^r$ and its limit~$\Phi_{P,L}^r$. The formula in the corollary can be compared with the descent identity (4.7) in~\cite{a-Ftr}, as suggested in the concluding remarks of~\cite{a-asym}. Note that the functions $\Phi_{M,L_1}^M$ have been determined in~\cite{he-1}.\medskip

{\noindent\it Example.} Suppose that $M=M_0$ is a minimal Levi subgroup. Then $L_1=M$, $\tau=\sigma\in\Pi_2(M)=\Pi\temp(M)$ is a finite-dimensional representation, and we need only consider $L\in\CL(M)$. By a theorem of Harish-Chandra's,
\[
\Phi_{M,M}^M(\gamma,\sigma)=|D^M(\sigma)|^{1/2}\tr\sigma(\gamma).
\]
Now $|D^M(\exp Y)|^{1/2}=\eps_\Sigma^U\Delta_\Sigma(\exp Y)$ for $Y$ in a connected component $U$ of $\exp^{-1}(T(\R)\cap M(\R)\reg)$, where $(\eps_\Sigma^U)^4=1$ and $\Delta_\Sigma$ denotes the denominator of the Weyl character formula with respect to a system $\Sigma$ of positive roots of $\t$ in~$\mathfrak m$. Thus, if $\mu\in\t^*(\sigma)$ is $\Sigma$-dominant, then
\[
\Phi_{M,M}^M(\exp Y,\sigma)=\eps_\Sigma^U\sum_{w\in W_T^M}\eps^M(w)e^{\mu(wY)}
\]
for $Y\in U$.
Note that we have a split exact sequence
\[
1\to W_T^M\to W_T^G\to W_M^G\to1,
\]
where the stabiliser of $\Sigma$ in $W_T^G$ is mapped bijectively onto~$W_M^G$ (cf.~\cite{ho-1}, Lemma~2). The formula in the corollary now specialises to
\begin{multline*}
\Phi_{P,L}(\exp Y,\sigma^L)
=n^L(\sigma)\eps_\Sigma^U
\sum_{w\in W_T^G}\eps^M(w)e^{\mu(wY)}\\
\times\sum_{S\in\CL(M)}d_M^G(L,S)m_M^S(\check\sigma,wP\cap S),
\end{multline*}
where $e^M(w)=(-1)^{\#(w\Sigma\cap-\Sigma)}$, provided $\sigma^L\in T\disc(L)$.

\section{The differential equations}

Let us recall from Proposition~11.1 of~\cite{a-loc} the differential equations satisfied by weighted orbital integrals. For every connected reductive $\R$-group~$G$, its Levi subgroup $M$ and maximal torus $T\subset M$, there exists a smooth map
\[
\partial_M=\partial_M^G:
T(\R)\cap G\reg(\R)\to\Hom(Z(\g),U(\t))
\]
such that, for all $\gamma\in T(\R)\cap G\reg(\R)$, $z\in Z(\g)$ and $f\in\CC(G)$,
\begin{equation}\label{diffeq}
J_M(\gamma,zf)=\sum_{S\in\CL(M)}\partial_M^S(\gamma,z_S)J_S(\gamma,f).
\end{equation}
Here we use a Haar measure on $G_\gamma(\R)=T(\R)$ independent of~$\gamma$. We denote by $z_S\in Z(\mathfrak s)$ the image of $z$ under the Harish-Chandra homomorphism, and the smoothness of the map $\partial_M(t,z)$ is meant for fixed $z$, where it takes values in a finite-dimensional space. The family of such maps is unique. Moreover, $\partial_G^G(z)=z_T$, which is the image of $z$ under the Harish-Chandra isomorphism $Z(\g)\to Z(\t)^W$, where $W$ is the Weyl group of~$(\g,\t)$.

The analogue of the differential equation~(\ref{diffeq}) with the weighted orbital integrals~$J_M(\gamma,f)$ replaced by the invariant distributions $I_M(\gamma,f)$ is equation~(2.6) of~\cite{a-invI}. As explained in~\cite{a-Ftr}, p.~197/198, it is also valid for $f\in\CC(G)$. We combine it with Theorem~4.1 of that paper, which we recalled in equation~(\ref{Ftr}). In that equation, we can interchange the differential operators with the integration due to the estimates (4.4) of~\cite{a-Ftr}. Moreover, $(z f)_L(\tau)=\check\chi(z)\hat f_L(\tau)$, where $\chi$ is the infinitesimal character of~$\tau^G$. We conclude that the functions $\Phi_M(\gamma)=\Phi_{M,L}(\gamma,\tau)$ satisfy the differential equations
\begin{equation}\label{diffeq1}
\chi(z)\Phi_M(\gamma)
=\sum_{S\in\CL(M)}\partial_M^{S}(\gamma,z_{S})\Phi_{S}(\gamma),
\end{equation}
because both left and right hand side satisfy the symmetry condition (4.2) of~\cite{a-Ftr}. If the equations are satisfied for $z\in\ker\chi$, then they are satisfied for all $z\in Z(\g)$, because $\chi(1)=1$ and $\partial_M^{S}(\gamma,1)=0$ for $M\ne S$.

The differential operators in~(\ref{diffeq1}) can be pulled back to the Lie algebra $\t(\R)$ of~$T(\R)$ under the exponential map and extended to meromorphic differential operators $\tilde\partial_M^S(Y,z_S)$ on its complexification $\t$.
For any parabolic $\R$-subgroup $P\supset T$ those differential operators are holomorphic on the sector
\[
\t_P=\{Y\in\t:|\alpha(Y)|>0\mbox{ for all roots $\alpha$ of $\t$ in~$\mathfrak n$}\},
\]
where $\mathfrak n$ denotes the (complex) Lie algebra of the unipotent radical $N$ of~$P$. Let $\CL(P)=\CL(M_P)$, where $M_P$ is the Levi component of $P$ that contains~$T$. A solution of the complexified system of differential equations
\begin{equation}\label{diffeq2}
\chi(z)\Psi_M(Y)
=\sum_{S\in\CL(M)}\tilde\partial_M^{S}(Y,z_{S})\Psi_{S}(Y)
\end{equation}
is then a tuple of holomorphic functions $\Psi_M$ on $\t_P$ indexed by the elements $M$ of~$\CL(P)$. The infinitesimal character $\chi$ corresponds to a $W$-orbit in $\t^*$ that we denote by~$\t^*(\chi)$. As a consequence of the theorem on the holomorphic dependence on parameters, the spaces of solutions $\CS_P(\chi)$ fit together to a holomorphic vector bundle $\CS_P$ over the affine space of characters of~$Z(\g)$, which we identify with $W\backslash\t^*$. For every subset $E$ of~$\t^*(\chi)$ we have the subset $\CS_P(\chi,E)$ of solutions whose leading exponents are in~$E$, and that is a subspace of $E$ is $P$-closed in~$\t^*(\chi)$ (i.~e., closed under subtraction of any sum, possibly with repetitions, of roots of $\t$ in $\mathfrak n$, cf.~\cite{ho-2}, p.~788). Let $\CE$ be an open subset of $\t^*$ such that its image under the natural map $\pi:\t^*\to W\backslash\t^*$ is also open and $\CE$ is closed in~$\pi^{-1}(\pi(\CE))$. If for each $\chi\in\pi(\CE)$ the set~$\t^*(\chi)\cap \CE$ is $P$-closed, then the spaces $\CS(\chi,\t^*(\chi)\cap \CE)$ make up a holomorphic subbundle $\CS_P(\CE)$ of~$\CS_P|_{\pi(\CE)}$.

We will now recall from Theorem~5.8 of~\cite{ho-2} the standard solutions. Suppose that $\mu$ is $P$-minimal (i.~e., $\{\mu\}$ is $P$-closed) in~$\t^*(\chi)$. Then for fixed $\gamma\in T(\R)^1$ there is a unique solution $(\Psi_M)_{M\in\CL(P)}$ of~(\ref{diffeq2}) such that
\begin{itemize}
\item $\Psi_G(Y)=e^{\mu(Y)}$,
\item If $M\ne G$, then there exist positive constants $C$ and $d$ such that
\[
\left|\Psi_M(Y)\right|\le C\|Y\|^d\max_\alpha\left|e^{(\mu-\alpha)(Y)}\right|,
\]
where $\alpha$ runs through the roots of~$\t$ in~$\mathfrak n$.
\end{itemize}
We refer to~\cite{ho-2} for the series expansion of this standard solution and explicit calculations for groups of low rank.

If the stabiliser $W_\mu$ of $\mu$ in~$W$ is nontrivial, there are further solutions. Indeed, if $c$ belongs to the space $H_{W_\mu}(\t^*)$ of $W_\mu$-harmonic elements of the symmetric algebra~$S(\t^*)$, then $(c\Psi_M)_{M\in\CL(P)}$ is also a solution, where $c$ is considered as a differential operator on $\t^*$ applied in the variable~$\mu$.

Since the standard solution $\Psi_M(Y)$ of~(\ref{diffeq2}) is uniquely determined by~$P$ and~$\mu$, we use the notation $\Psi_M^P(Y,\mu)$. This function is holomorphic in both arguments as the limit of a normally convergent series of holomorphic functions.
Given $L\in\CL(P)$, we obtain another solution of the system~(\ref{diffeq2}) by restricting $\Psi_M^{P\cap L}$ to $\t_P$ for $M\subset L$ and setting the components with $M\not\subset L$ equal to zero.

According to Theorem~5.8 of~\cite{ho-2}, the standard solutions $\Psi_M^P$ are the building blocks of arbitrary solutions of the system of differential equations~(\ref{diffeq1}) in the following sense. If every element of $\t^*(\chi)$ is $P$-minimal in this set, then for any solution $(\Phi_M)_{M\in\CL(P)}$ and any connected open subset $U$ of $\exp^{-1}(T(\R)\cap P\reg(\R))$
there are unique elements $c_M^{P,U}(\mu)\in H_{W_\mu}(\t^*)$ for every $M\in\CL(P)$ such that
\begin{equation}\label{defc0}
\Phi_M(\exp Y)=\sum_{\mu\in\t^*(\chi)}\sum_{S\in\CL(M)}
c_{S}^{P,U}(\mu)\Psi_M^{P\cap S}(Y,\mu)
\end{equation}
for all $M\in\CL(P)$ and all $Y\in U_P=U\cap\t_P$. Note that $\exp U_P\subset G\reg(\R)$. Actually, the proof was only given for regular infinitesimal characters but extends easily to this slightly more general situation, as for dimensional reasons the standard solutions, augmented by harmonic differential operators as above, still span the full space of solutions. We will write $\CU^P$ for the set of connected components of~$\exp^{-1}(T(\R)\cap M_{P,\mathrm{reg}}(\R))$.

There is a simple version of parabolic descent. If $Q$ is a parabolic subgroup containing~$P$, then $\t_P\subset\t_Q$, and for any Levi subgroup $M$ containing the Levi component of~$Q$ that contains $T$ we conclude from the characterisation of $\Psi_M^P$ that
\[
\Psi_M^Q\big|_{\t_P}=\Psi_M^P.
\]
For $\mu\in\t^*(\chi)$ that is $P$-minimal in this set we deduce the equality
\[
c_M^{P,U}(\mu)=c_M^{Q,U}(\mu)
\]
for $U\in\CU^Q$ using the uniqueness of these numbers. The applications of this equality are limited. Given $M\in\CL(P)$, we cannot always find $Q\in\CP(M)$ such that $P\subset Q$. We defer a more detailed study of parabolic descent to the next section.

Now we will apply the quoted results to the Fourier transforms $\Phi_{M,L}(\gamma,\tau)$. A virtual character $\tau\in T\disc(L)$ will be called $P$-semiregular if every element of $\t^*(\tau^G)$ is $P$-minimal in that set, and we write $T\disc^{P-\mathrm{sreg}}(L)$ for the corresponding open subset. We see that for every $M\in\CL(P)$, $U\in\CU^P$, $\tau\in T\disc^{P-\mathrm{sreg}}(L)$ and $\mu\in\t^*(\tau^G)$ there are unique elements $c_{M,L}^{P,U}(\tau,\mu)\in H_{W_\mu}(\t^*)$ such that
\begin{equation}\label{defc}
\Phi_{M,L}(\exp Y,\tau)=\sum_{\mu\in\t^*(\tau^G)}\sum_{S\in\CL(M)}
c_{S,L}^{P,U}(\tau,\mu)\Psi_M^{P\cap S}(Y,\mu)
\end{equation}
for $Y\in U_P$.

An element $\tau\in T\disc^{P-\mathrm{sreg}}(L)$ will be called $P$-regular if $|W_\mu|$ is minimal for every $\mu\in\t^*(\tau^G)$ and $\tau^{L'}\notin T\disc(L')$ for any Levi subgroup $L'$ properly containing~$L$. Such elements make up a dense open subset $T\disc^{P-\mathrm{reg}}(L)$, over which we have a smooth covering $\{(\tau,\mu)\mid\mu\in\t^*(\tau^G)\}$, and the functions $c_{M,L}^{P,U}(\tau,\mu)$ are smooth on this covering.

\begin{theo}\label{exponents}
If $M\in\CL(P)$, $U\in\CU^P$, $\tau\in T\disc^{P-\mathrm{sreg}}(L)$ and $\mu\in\t^*(\tau^G)$, then $c_{M,L}^{P,U}(\tau,\mu)$ vanishes unless $\Re\mu(X)\le0$ for all $X\in\a_P^+$.
\end{theo}
\begin{proof}
We may assume by induction that the theorem is true for all Levi subgroups $S$ properly containing~$M$ and for all $\tau'\in T\disc(L')$ such that $L'$ is properly contained in~$L$ and $\tau'^L=\tau$. The basis of induction is a special case of the inductive step as for $M=G$ there are just no such~$S$ and for $\tau\in T\ell(L)$ there are no such~$(L',\tau')$.

We fix $U$, $L$ and $\tau\in T\disc^{P-\mathrm{reg}}(L)$. Let $\phi$ be a smooth function on~$T\disc(L)$ not vanishing at~$\tau$ whose support is compact and contained in the $i\a_L^*$-orbit of~$\tau$ intersected with $T\disc^{P-\mathrm{reg}}(L)$. There exists a Levi subgroup $L_1\subset L$ such that $\tau$ is induced from $T\ell(L_1)$. If the support of~$\phi$ is small enough, then due to the regularity of~$\tau$ the pull-back of $\phi$ under induction extends to a smooth compactly supported function $\phi_1$ on $T\ell(L_1)$ that vanishes on the preimage of $T\disc(L')$ for $L'\in\CL(L_1)$ unless $L'\subset L$. Due to the trace Paley-Wiener theorem (see~\cite{a-PW}), there is a function $f\in\CC(G)$ such that $\hat f_{L'}=0$ on $T\ell(L')$ unless $L'$ is conjugate to~$L_1$, while
\[
\hat f_{L_1}(\tau_1)=\frac1{|W_{L_1}|}\sum_{w\in W_{L_1}}\phi_1(w\check\tau_1)
\]
for $\tau_1\in T\ell(L_1)$. By construction, $\hat f_{L'}=0$ on $T\disc(L')$ unless $L'$ is conjugate to a Levi subgroup sandwiched between $L_1$ and~$L$. According to~(\ref{Ftr}) we have
\[
I_M(\exp Y,f)=\sum_{[L']}
\int_{W_{L'}\backslash T\disc(L')}
\Phi_{M,L'}(\exp Y,\tau')\hat f_{L'}(\check\tau')\,d\tau'
\]
for $Y\in U_P$, where the sum is taken over all $G(\R)$-conjugacy classes of Levi subgroups which have a representative $L'$ sandwiched between $L_1$ and~$L$. If we replace $Y$ by $Y+X$, where $X\in\a_P^+=\a_M\cap\t_P$,  the left-hand side is rapidly decreasing as $X\to\infty$ by Corollary~7.4 of~\cite{a-ds} or rather its analogue for the invariant distribution $I_M$ (cf.~\cite{a-Ftr}, p.~213). Due to the inductive assumption, the same is true for the terms on the right-hand side with $L'\ne L$ and for the contribution of Levi subgroups $S\ne M$ to~(\ref{defc}) in view of the definition of~$\Psi_M^P$ and the continuity of $c_{S,L}^P$ on the support of~$\phi$. We conclude that
\[
g(Y)=\int_{i\a_L^*}
\sum_{\mu\in\t^*(\tau_\lambda^G)}
c_{M,L}^{P,U}(\tau_\lambda,\mu,Y)e^{\mu(Y)}
\phi(\tau_\lambda)\,d\lambda
\]
is rapidly decreasing as well. Here we have interpreted $c_{M,L}^{P,U}(\tau_\lambda,\mu)$ as a harmonic polynomial on~$\t$ and inserted the additional argument~$Y$.

Let us fix a maximal torus $T_1$ of $L$ and a parameter $\mu_1\in\t_1^*(\tau)$. We denote by $V$ the set of isomorphisms $\t\to\t_1$ induced by elements of~$G(\C)$. Then for $\lambda\in i\a_L^*$ we have
\[
\t^*(\tau_\lambda^G)=\{v^*(\mu_1+\lambda)\mid v\in V\}.
\]
We fix $X$ and $Y$ as above. With the notation
\[
c_v(\lambda,t)=e^{(\mu_1+\lambda)(vY)}
c_{M,L}^{P,U}(\tau_\lambda,v^*(\mu_1+\lambda),tX+Y)
\]
we get
\begin{align*}
g(tX+Y)&=\sum_{v\in V}\int_{i\a_L^*}
c_v(\lambda,t)e^{t(\mu_1+\lambda)(vX)}\phi(\tau_\lambda)\,d\lambda\\
&=\sum_{v\in V}\int_{i\a_L^*}
\phi(\tau_\lambda)
c_v(\lambda)e^{at}\big|_{a=(\mu_1+\lambda)(vX)}\,d\lambda,
\end{align*}
where the polynomial $c_v(\lambda)$ in the variable $t$ has been reinterpreted as a differential operator acting in the variable~$a$. The left-hand side is rapidly decreasing for $t\to\infty$, so the Laplace transform
\[
G(s)=\int_0^\infty g(tX+Y)e^{-st}\,dt
\]
is holomorphic for $\Re s>0$. If we plug in the formula for $g$, we can interchange the order of integration for $\Re s\gg0$ and get
\[
G(s)=\sum_{v\in V}\int_{i\a_L^*}
\phi(\tau_\lambda)
\left.c_v(\lambda)\frac1{s-a}\right|_{a=(\mu_1+\lambda)(vX)}\,d\lambda.
\]
In our case the inversion formula
\[
g(tX+Y)=\frac1{2\pi i}\lim_{R\to\infty}
\int_{\sigma-iR}^{\sigma+iR}G(s)e^{st}\,ds
\]
is valid for $t>0$ and $\sigma>0$.

Let $V_+=\{v\in V\mid\Re\mu_1(vX)>0\}$ and $\eps=\min\{\Re\mu_1(vX)\mid v\in V_+\}$. We may assume that the support of $\phi$ is so small that for any of its elements $\lambda$ and any $v\in V$ we have $|\lambda(vX)|<\eps/2$. Then the term in the formula for $G$ corresponding to a given $v$ extends to a holomorphic function outside the disc about $\mu_1(vX)$ with radius~$\eps/2$. We may now interchange summation and integration in the inversion formula if $\sigma\gg0$ or if $\sigma=\eps/2$. Since the result is independent of~$\sigma$, we conclude that
\[
\int_{i\a_L^*}\sum_{v\in V_+}
c_v(\lambda,t)e^{t(\mu_1+\lambda)(vX)}\phi(\tau_\lambda)\,d\lambda=0.
\]
The function $\phi$ was arbitrary in a neighbourhood of~$\tau$, thus
\[
\sum_{v\in V_+}
c_{M,L}^{P,U}(\tau,v^*\mu_1,Y)
e^{\mu_1(v(Y))}=0
\]
for any $Y\in U_P$. Since $c$ depends polynomially on~$Y$, the terms of the sum are linearly independent as functions on~$U_P$. Consequently, each of them must vanish, and the theorem is proved in the case $\tau\in T\disc^{P-\mathrm{reg}}(L)$.

If $\CE$ is an open neighbourhood of $\{\mu\in\t^*\,:\,\Re\mu(X)\le0\mbox{ for all $X\in\a_P^+$}\}$, then it follows that for all elements of $T\disc^{P-\mathrm{reg}}(L)$ the tuple of functions $\Phi_{M,L}(U_X,\tau)$ lies in the subbundle~$\CS_P(\CE)$ (or vanishes if $\t^*(\tau)\notin\pi(\CE)$). By continuity, this is also the case for all $\tau\in T\disc^{P-\mathrm{sreg}}(L)$.\end{proof}

Now we are going to apply the asymptotic formula to the expression~(\ref{defc}).
\begin{theo}\label{asymPhi}
Let $T$ be a maximal torus in the Levi component~$M$ of the parabolic subgroup $P$ of $G$. Let $\tau\in T\ell^{P-\mathrm{reg}}(L)$ be represented by $(M_1,\sigma,r)$, where $\sigma\in\Pi_2(M_1)$ and $r\in R_\sigma^L$. For $U\in\CU^P$ and $Y\in U_P$, we have
\[
\Phi_{P,L}(\exp Y,\tau)
=\frac{|W_M|}{|W_{M_1}|}
\sum_{\mu\in\t^*(\tau^G,M)}
c_{M,L}^{P,U}(\tau,\mu)e^{\mu(Y)},
\]
where $\t^*(\tau^G,M)$ denotes the set of all $\mu\in\t^*(\tau^G)$ for which there exists a $G(\R)$-conjugate $(M_1',\sigma')$ of $(M,\sigma)$ such that $M_1'\subset M$ and $\mu|_{\a_M}$ is the infinitesimal central character of~$\sigma'^M$, and for which this remains true in a neighbourhood of~$\tau$.
\end{theo}
\begin{proof}
As in the previous proof, we fix $L$, $\tau$ and $U$ and a function $\phi$, which gives rise to a function $f\in\CC(G)$. If we plug in equation~(\ref{defc}) into equation~(\ref{Ftr}) with $(\gamma,f)$ replaced by~$(\exp Z,f_X)$, where $Z=X+Y$, $X\in\a_P^+$ and $Y\in U_P$, we get
\begin{multline*}
I_M(\exp Z,f_X)=
\int_{i\a_L^*}\sum_{\mu\in\t^*(\tau_\lambda^G)}
\sum_{S\in\CL(M)}\!\!
c_{S,L}^{P,U}(\tau,\mu)\Psi_M^{P\cap S}(\gamma,Z,\mu)
\hat\alpha_X(\check\tau_\lambda)\phi(\tau_\lambda)\,d\lambda.
\end{multline*}
We know from Theorem~\ref{exponents} that the coefficients $c_{S,L}^{P,U}(\tau,\mu)$ vanish unless $e^{\mu(X)}$ is bounded for~$X\in\a_P^+$. In this case, the functions $\Psi_M^{P\cap S}(Z,\mu)$ with $S\ne M$ are exponentially decreasing as $X\mathrel{\mathop{\to}\limits_P}\infty$. Moreover, $|\hat\alpha_X(\check\tau)|\le1$ due to the fact that the central character of a unitary representation is unitary. Since the functions $c_{S,L}^P$ are smooth on the support of~$\phi$, we can take the sum over $S$ out of the integral, and the terms with $S\ne M$ tend to zero by bounded convergence. Thus we are left with
\begin{multline*}
\lim_{X\mathrel{\mathop{\to}\limits_{P,r}}\infty}
I_M(\exp Z,f_X)\\=
\lim_{X\mathrel{\mathop{\to}\limits_{P,r}}\infty}
\int_{i\a_L^*}\sum_{\mu\in\t^*(\tau_\lambda^G)}
c_{M,L}^{P,U}(\tau_\lambda,\mu,Z)e^{\mu(Z)}
\hat\alpha_X(\check\tau_\lambda)\phi(\tau_\lambda)\,d\lambda,
\end{multline*}
where we have interpreted $c_{M,L}^P$ as a harmonic polynomial in the variable~$Z$ again. Now we replace $X$ by $tX$ with fixed $X\in\a_P^+$ and $Y\in U_P$. We rewrite the left-hand side in the notation of Corollary~\ref{PhiP} and the right-hand side in the notation from the proof of Theorem~\ref{exponents}:
\begin{multline*}
\int_{i\a_L^*}\Phi_{P,L}(\exp Y,\tau_\lambda)
\phi(\tau_\lambda)\,d\lambda\\
=\lim_{t\to\infty}
\sum_{v\in V}\int_{i\a_L^*}
c_v(\lambda,t)e^{t(\mu_1+\lambda)(vX)}
\hat\alpha_X(\check\tau_\lambda)\phi(\tau_\lambda)\,d\lambda.
\end{multline*}

If $\Re(\mu_1+\lambda)(vX)<0$ for some $X\in\a_P^+$ and some inner point $\tau_\lambda$ of the support of~$\phi$ then the contribution of $v$ tends to zero as $t\to\infty$ by bounded convergence. Therefore we are left with the terms for those $v$ only in which $\Re\mu(vX)=0$ and $\Re\lambda(vX)=0$ for all $\lambda\in i\a_L^*$. Choosing a generic~$X$, we see that only those $v$ contribute to the limit for which $\mu_1$ and all $\lambda\in i\a_L^*$ have purely imaginary values on~$v\a_M$. Let us denote the set of such $v$ by~$V_1$.

Let $\nu\in i\a_{M_1}^*$ be the infinitesimal central character of $\sigma\in\Pi_2(M_1)$ as in the statement of the theorem. Then we have
\[
\hat\alpha_X(\check\tau_\lambda)=\frac1{|W_{M_1}|}\sum_{u\in U(M,M_1)}e^{-(\nu+\lambda)(uX)},
\]
where $U(M,M_1)$ denotes the set of embeddings $\a_M\to\a_{M_1}$ induced by elements of~$G(\R)$ as in~\cite{a-asym},~\S1. Since $X$ is in the span of the cocharacters of~$T$, its image under $v$ is in the span of the cocharacters of~$T_1$.  We may assume that $T_1\subset M_1$, so that $\Im(\mu_1+\lambda)(vX)=(\nu+\lambda)(vX)$ for all $\lambda\in i\a_L^*$ and
\begin{multline*}
\int_{i\a_L^*}\Phi_{P,L}(\exp Y,\tau_\lambda)
\phi(\tau_\lambda)\,d\lambda\\
=\lim_{t\to\infty}
\sum_{v\in V_1}\frac1{|W_{M_1}|}\sum_{u\in U(M,M_1)}e^{t\nu(vX-uX)}
\int_{i\a_L^*}
c_v(\lambda,t)e^{t\lambda(vX-uX)}
\phi(\tau_\lambda)\,d\lambda.
\end{multline*}
Since $c_v\phi$ is smooth and compactly supported in~$\lambda$ and depends polynomially on~$t$, the term for given $u$ and~$v$ tends to zero unless $vX-uX$ vanishes on~$\a_L^*$. Choosing $X$ generic, we see that this condition has to be satisfied for all~$X\in\a_M$. Let $V_2$ be the set of all $v\in V_1$ for which there exists $u\in U(M,M_1)$ so that it is satisfied. It will then be satisfied for the whole coset $uW_M$, which is all of~$U(M,M_1)$. By Lemma~\ref{Tdisc}, we have $\nu\in i\a_L^*$, so that the exponential term outside the integral disappears as well, and we get in the original notation
\begin{multline*}
\int_{i\a_L^*}\Phi_{P,L}(\exp Y,\tau_\lambda)
\phi(\tau_\lambda)\,d\lambda\\
=\lim_{t\to\infty}\sum_{v\in V_2}\frac{|W_M|}{|W_{M_1}|}
\int_{i\a_L^*}
c_{M,L}^{P,U}(\tau_\lambda,v^*(\mu_1+\lambda),tX+Y)e^{(\mu_1+\lambda)(vY)}
\phi(\tau_\lambda)\,d\lambda.
\end{multline*}
The right-hand side is a polynomial in~$t$, and from the existence of the limit it follows that it does not depend on~$t$. The terms for various $v$ are linearly independent as functions of~$Y$, so each of them must be independent of~$t$. Since $\phi$ was an arbitrary smooth function of compact support in a neighbourhood of~$\tau$, the asserted equality follows. The definition of $V_2$ can be rephrased as the condition on~$\mu$ in the definition of~$\t^*(\tau^G,M)$.
\end{proof}

If we compare the expressions for $\Phi_{P,L}(\exp Y,\sigma)$ given in Corollary~\ref{PhiP} and Theorem~\ref{asymPhi}, we can determine part of the coefficients~$c_{M,L}^{P,U}$.\medskip

{\noindent\it Example.} In the example given at the end of the preceding section, let us specialise to the case that $\sigma^L\in T\ell(L)$. Then we have for all $w\in W_T^G$ that $w\mu\in\t^*(\sigma,M)$ and
\[
c_{M,L}^{P,U}(\sigma^L,w\mu)=n^L(\sigma)\eps_\Sigma^U\eps^M(w)
\sum_{S\in\CL(M)}d_M^G(L,S)m_M^S(\check\sigma,w^{-1}P\cap S).
\]
In particular,
\[
c_{M,M}^{P,U}(\sigma,w\mu)=\eps_\Sigma^U\eps^M(w)m_M(w\check\sigma,P).
\]

\section{Parabolic descent}

Weighted orbital integrals satisfy descent identities. If $M\subset L$ are Levi subgroups of $G$ and $\gamma\in M(\R)\cap G\reg(\R)$, then according to equation~(1.5) of~\cite{a-Ftr} we have
\begin{equation}\label{decJ}
J_L(\gamma,f)=\sum_{S\in\CL(M)}d_M^G(L,S)J_M^S(\gamma,f_{Q_S}),
\end{equation}
where we can choose any $Q_S\in\CP(S)$ for each~$S$. We are going to prove parallel descent identities for the differential equations satisfied by the distributions $I_M$ and their Fourier transforms, for the standard solutions of those differential equations and for the coefficients in the expression of the Fourier transforms in terms of standard solutions.

\begin{prop}\label{decdiff}
  Let $M\subset L$ be Levi subgroups of $G$ and $T\subset M$ a maximal torus. Then, for all $\gamma\in T(\R)\cap G\reg(\R)$ and $z\in Z(\g)$, we have
\[
\partial_L(\gamma,z)=\sum_{S\in\CL(M)}d_M^G(L,S)\partial_M^S(\gamma,z_S).
\]
\end{prop}
\begin{proof}
   We change variables for Levi subgroups in favour of a more systematic notation later in the proof and thus apply equation~(\ref{decJ}) to Levi subgroups $M_1\subset M$ containing~$T$, viz.
\[
J_M(\gamma,f)=\sum_{G_1\in\CL(M_1)}d_{M_1}^G(M,G_1)
J_{M_1}^{G_1}(\gamma,f_{Q_{G_1}}).
\]
Replacing $f$ by~$zf$ and observing that $(zf)_{Q_{G_1}}=z_{G_1}f_{Q_{G_1}}$, we can apply equation~(\ref{diffeq}) to the right-hand side, and with the transitivity of the Harish-Chandra homomorphism we get
\[
J_M(\gamma,zf)=\sum_{G_1\in\CL(M_1)}d_{M_1}^G(M,G_1)
\sum_{S_1\in\CL^{G_1}(M_1)}\partial_{M_1}^{S_1}(\gamma,z_{S_1})
J_{S_1}^{G_1}(\gamma,f_{Q_{G_1}}).
\]
If $d_{M_1}^G(M,G_1)\ne0$, then $\a_{M_1}^M\cap\a_{M_1}^{G_1}=0$ and hence $\a_{M_1}^M\cap\a_{M_1}^{S_1}=0$ for $S_1\in\CL^{G_1}(M_1)$. Given such~$S_1$ and any $S\in\CL(M)\cap\CL(S_1)$, we have $d_{M_1}^S(M,S_1)\ne0$ if and only if $S$ is the unique Levi subgroup such that $\a_{M_1}^S=\a_{M_1}^M\oplus\a_{M_1}^{S_1}$, and thus
\begin{equation}\label{trand}
d_{M_1}^G(M,G_1)=\sum_{S\in\CL(M)\cap\CL(S_1)}d_{M_1}^S(M,S_1)d_{S_1}^G(S,G_1).
\end{equation}
The Hasse diagram may help keeping track of the various inclusions.
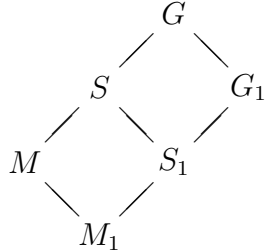
\begin{center}
\setlength{\unitlength}{1.2em}
\newcommand{\ul}{\makebox(1,1){\line(1,1){0.9}}}
\newcommand{\dl}{\makebox(1,1){\line(1,-1){0.9}}}
\begin{picture}(7,7)
\put(0,2){\makebox(1,1){$M$}}
\put(2,4){\makebox(1,1){$S$}}
\put(4,6){\makebox(1,1){$G$}}
\put(2,0){\makebox(1,1){$M_1$}}
\put(4,2){\makebox(1,1){$S_1$}}
\put(6,4){\makebox(1,1){$G_1$}}
\put(1,1){\dl} \put(3,3){\dl} \put(5,5){\dl}
\put(1,3){\ul} \put(3,5){\ul}
\put(3,1){\ul} \put(5,3){\ul}
\end{picture}
\end{center}
After changing the order of summation, the expression for $J_M(\gamma,zf)$ becomes
\[
\sum_{S\in\CL(M)}\sum_{S_1\in\CL^S(M_1)}d_{M_1}^S(M,S_1)
\partial_{M_1}^{S_1}(\gamma,z_{S_1})
\sum_{G_1\in\CL(S_1)}d_{S_1}^G(S,G_1)J_{S_1}^{G_1}(\gamma,f_{Q_{G_1}}).
\]
Using the descent identity for weighted orbital integrals again, we get
\[
J_M(\gamma,zf)=\sum_{S\in\CL(M)}\sum_{S_1\in\CL^S(M_1)}d_{M_1}^S(M,S_1)
\partial_{M_1}^{S_1}(\gamma,z_{S_1})J_S(\gamma,f).
\]
Since the Theorem is trivially true for $\R$-anisotropic groups~$G$, we may inductively assume that it is true for $G$ replaced by a proper Levi subgroup. If we subtract the last equality from equation~(\ref{diffeq}), we see that all the terms with $S\ne G$ cancel, and we are left with
\[
\left(\partial_M(\gamma,z)-\sum_{S_1\in\CL(M_1)}d_{M_1}(M,S_1)
\partial_{M_1}^{S_1}(\gamma,z_{S_1})\right)J_G(\gamma,f)=0.
\]
If $\gamma_0\in T(\R)\cap G\reg(\R)$, then the stabiliser of $\gamma_0$ in $W_T$ is trivial. Thus any smooth function in a sufficiently small neighbourhood of $\gamma_0$ in~$T$ is of the form $J_G(\gamma,f)$ for a function $f\in\CH(G)$, and our assertion follows.
\end{proof}

Next we consider the standard solutions of the system of differential equations~(\ref{diffeq2}). We fix a parabolic subgroup $P$ of $G$ and a maximal torus $T\subset P$ and denote again by $\CL(P)$ the set of Levi subgroups of $G$ containing the Levi component of $P$ that contains~$T$.
\begin{prop}\label{decPsi}
For all Levi subgroups $M\subset L$ in~$\CL(P)$,
$Y\in\t_P$ and any $\mu\in\t^*$ that is $P$-minimal in its $W$-orbit, we have
\[
\Psi_L^P(Y,\mu)=\sum_{S\in\CL(M)}d_M^G(L,S)\Psi_M^{P\cap S}(Y,\mu).
\]
\end{prop}
\begin{proof}
Since the standard solutions depend holomorphically on~$\mu$ subject to the assumption of the Theorem, it suffices to consider regular $\mu$ which take purely imaginary values on~$\t(\R)$. As in the preceding proof, we will change the notation $(M,L)$ to~$(M_1,M)$. We may assume inductively that the Theorem is true for $M$ replaced by any Levi subgroup $L$ of~$G$ properly containing~$M$. The basis of induction is a special case of the inductive step, as for $M=G$ there are just no such~$L$.

Let $z\in\ker\chi$, where $\chi$ is the infinitesimal character parametrised by~$\mu$. Then, for $Y\in\t_P$,
\[
\sum_{S\in\CL(M)}\tilde\partial_M^S(Y,z_S)
\Psi_S^P(Y,\mu)=0.
\]
Let $M_1\in\CL(P)$ be a Levi subgroup of $M$. From Theorem~\ref{decdiff} we get, after pull-back under the exponential map and meromorphic extension,
\[
\tilde\partial_M^S(Y,z_S)=\sum_{S_1\in\CL^S(M_1)}d_{M_1}^S(M,S_1)
\tilde\partial_{M_1}^{S_1}(Y,z_{S_1})=0.
\]
Starting with the pair $M\subset G$ replaced by Levi subgroups $M_1\subset G_1$, we obtain for all $Y\in\t_{P\cap G_1}$
\[
\sum_{S_1\in\CL^{G_1}(M_1)}
\tilde\partial_{M_1}^{S_1}(Y,z_{S_1})
\Psi_{S_1}^{P\cap G_1}(Y,\mu).
\]
We restrict $Y$ to the set $\t_P$, multiply the equation with $d_{M_1}^G(M,G_1)$ and sum over all $G_1\in\CL(M_1)$. Using equation~(\ref{trand}), we can introduce an additional summation over $S\in\CL(M)\cap\CL(S_1)$ just as in the proof of Theorem~\ref{decdiff}. After changing the order of summation, we obtain the vanishing of
\[
\sum_{S\in\CL(M)}\sum_{S_1\in\CL^S(M_1)}d_{M_1}^S(M,S_1)
\tilde\partial_{M_1}^{S_1}(Y,z_{S_1})
\sum_{G_1\in\CL(S_1)}d_{S_1}^G(S,G_1)
\Psi_{S_1}^{P\cap G_1}(Y,\mu).
\]
Forming the difference with our earlier result, we get
\begin{multline*}
\sum_{S\in\CL(M)}\sum_{S_1\in\CL^S(M_1)}d_{M_1}^S(M,S_1)\\
\times\tilde\partial_{M_1}^{S_1}(Y,z_{S_1})
\left(\Psi_S^P(Y,\mu)-\sum_{G_1\in\CL(S_1)}d_{S_1}^G(S,G_1)
\Psi_{S_1}^{P\cap G_1}(Y,\mu)\right)=0.
\end{multline*}
Due to the inductive assumption, all terms with $S\ne M$ vanish. But when $S=M$ and $d_{M_1}^S(M,S_1)\ne0$, then $S_1=M_1$, and we are left with
\[
z_T\left(\Psi_M^P(Y,\mu)-\sum_{G_1\in\CL(M_1)}d_{M_1}^G(M,G_1)
\Psi_{M_1}^{P\cap G_1}(Y,\mu)\right)=0.
\]
A function on $\t_P$ annihilated by all $z_T$ for $z\in\ker\chi$ is a linear combination of functions of the form $e^{w\mu(Y)}$ with $w\in W$ if $\mu$ is regular. Under the further assumption on $\mu$ made at the beginning of the proof, these functions have absolute value~1 for $Y\in\t(\R)$. For $M\ne G$, the expression in the brackets multiplied by $e^{-\mu(Y)}$ tends to zero as $Y\mathrel{\mathop{\to}\limits_P}\infty$ by definition of the standard solutions, so it must vanish identically. For $M=G$ this is also true by definition.
\end{proof}

Finally, we turn to the coefficients in the expression~(\ref{defc}) of Fourier transforms in terms of standard solutions. We keep $P$ and $T$ fixed.
\begin{prop}\label{decc}
Let $U\in\CU^P$ and $\tau\in T\disc^{P-\mathrm{sreg}}(L)$, where $L$ is a Levi subgroup of~$G$. Then, for any Levi subgroups $M_1\subset M$ in~$\CL(P)$, we have
\[
c_{M,L}^{P,U}(\tau,\mu)
=\sum_{S\in\CL(M_1)}d_{M_1}^G(M,S)
\sum_{\substack{w\in W_0^S\backslash W_0^G\\wL\subset S\\\mu\in\t^*(w\tau^S)}}c_{M_1,wL}^{P\cap S,U}(w\tau,\mu).
\]
\end{prop}
\begin{proof}
  On one hand, we can plug in the descent formula from Theorem~\ref{decPsi} into equation~(\ref{defc}) and get for $Y\in U_P$
\begin{multline*}
\Phi_{M,L}(\exp Y,\tau)=\sum_{\mu\in\t^*(\tau^G)}\sum_{S\in\CL(M)}
\sum_{S_1\in\CL^{S}(M_1)}d_{M_1}^{S}(M,S_1)\\
\times c_{S,L}^{P,U}(\tau,\mu)\Psi_{M_1}^{P\cap S_1}(Y,\mu).
\end{multline*}
On the other hand, we can use the descent formula (4.3) of~\cite{a-Ftr} for the Fourier transforms of the invariant distributions $I_M$, which reads in our notation
\[
\Phi_{M,L}(\gamma,\tau)=\sum_{G_1\in\CL(M_1)}d_{M_1}^G(M,G_1)\sum_{\substack{w\in W_0^{G_1}\backslash W_0^G\\wL\subset G_1}}\Phi_{M_1,wL}^{G_1}(\gamma,w\tau).
\]
Here we plug in the analogue of equation~(\ref{defc}), namely
\[
\Phi_{M_1,wL}^{G_1}(\exp Y,w\tau)
=\sum_{\mu\in\t^*(w\tau^{G_1})}\sum_{S_1\in\CL^{G_1}(M_1)}
c_{S_1,wL}^{P\cap G_1,U}(w\tau,\mu)\Psi_{M_1}^{P\cap S_1}(Y,\mu).
\]
If we change the order of summation, we obtain
\begin{multline*}
\Phi_{M,L}(\exp Y,\tau)=\sum_{\mu\in\t^*(\tau^G)}\sum_{S_1\in\CL(M_1)}
\sum_{G_1\in\CL(S_1)}d_{M_1}^G(M,G_1)\\
\sum_{\substack{w\in W_0^{G_1}\backslash W_0^G\\wL\subset G_1\\\mu\in\t^*(w\tau^{G_1})}}
c_{S_1,wL}^{P\cap G_1,U}(w\tau,\mu)
\Psi_{M_1}^{P\cap S_1}(Y,\mu).
\end{multline*}
Next we plug in equation~(\ref{trand}) and change the order of summation once more. If we compare the two expressions for $\Phi_{M,L}(\exp Y,\tau)$, we see that
\begin{multline*}
  \sum_{\mu\in\t^*(\tau^G)}\sum_{S\in\CL(M)}
\sum_{S_1\in\CL^{S}(M_1)}d_{M_1}^{S}(M,S_1)\\
\times\left(c_{S,L}^{P,U}(\tau,\mu)-\sum_{G_1\in\CL(S_1)}d_{S_1}^G(S,G_1)
\sum_{\substack{w\in W_0^{G_1}\backslash W_0^G\\
wL\subset G_1\\\mu\in\t^*(w\tau^{G_1})}}
c_{S_1,wL}^{P\cap G_1,U}(w\tau,\mu)\right)\\
\times\Psi_{M_1}^{P\cap S_1}(Y,\mu)=0.
\end{multline*}
We may assume inductively that the Theorem is true with $M$ replaced by a Levi subgroup $S$ properly containing~$M$. The basis of induction is a special case of the inductive step, as for $M=G$ there are just no such~$S$. Using the inductive assumption, we are left with the terms in which $S=M$ and hence $S_1=M_1$, i.~e.
\begin{multline*}
\sum_{\mu\in\t^*(\tau^G)}\left(c_{M,L}^{P,U}(\tau,\mu)
-\sum_{G_1\in\CL(M_1)}d_{M_1}^G(M,G_1)
\sum_{\substack{w\in W_0^{G_1}\backslash W_0^G\\wL\subset G_1\\\mu\in\t^*(w\tau^{G_1})}}
c_{M_1,wL}^{P\cap G_1,U}(w\tau,\mu)\right)\\
\times e^{\mu(Y)}=0.
\end{multline*}
We may again interpret the coefficients $c$ as polynomials of the variable~$Y$ rather than differential operators on~$\t_\C^*$. Exponential functions with different exponents, multiplied by polynomials, are linearly independent on~$U_P$.
\end{proof}

Finally, we can deduce a version of the descent identity~(\ref{defc}) in which the roles of $c$ and $\Psi$ have been switched.
\begin{cor}
For $Y\in U_P$ and $\tau\in T\disc^{P-\mathrm{reg}}(L)$
as in Theorem~\ref{decc}, we have
  \[
\Phi_{M,L}(\exp Y,\tau)=\sum_{S\in\CL(M)}
\sum_{\substack{w\in W_0^S\backslash W_0^G\\wL\subset S}}
\sum_{\mu\in\t^*(w\tau^S)}
c_{M,wL}^{P\cap S,U}(w\tau,\mu)\Psi_S^P(Y,\mu).
\]
\end{cor}
To see this, we rewrite equation~(\ref{defc}) for a Levi subgroups $S_1$ in place of~$M$ as
\[
\Phi_{S_1,L}(\exp Y,\tau)
=\sum_{\mu\in\t^*(\tau^G)}\sum_{G_1\in\CL(S_1)}
c_{G_1,L}^{P,U}(\tau,\mu)\Psi_{S_1}^{P\cap G_1}(Y,\mu)
\]
and plug in the descent identity from Theorem~\ref{decc} in the form
\[
c_{G_1,L}^{P,U}(\tau,\mu)
=\sum_{M\in\CL(M_1)}d_{M_1}^G(M,G_1)
\sum_{\substack{w\in W_0^M\backslash W_0^G\\wL\subset M\\\mu\in\t^*(w\tau^M)}}c_{M_1,wL}^{P\cap M,U}(w\tau,\mu).
\]
Here we plug in equation~(\ref{trand}) and change the order of summation to obtain
\begin{multline*}
\Phi_{S_1,L}(\exp Y,\tau)
=\sum_{\mu\in\t^*(\tau^G)}\sum_{S\in\CL(S_1)}\\
\sum_{M\in\CL^S(M_1)}d_{M_1}^S(M,S_1)
\sum_{\substack{w\in W_0^M\backslash W_0^G\\
wL\subset M\\\mu\in\t^*(w\tau^M)}}
c_{M_1,wL}^{P\cap M,U}(w\tau,\mu)\\
\times\sum_{G_1\in\CL(S_1)}d_{S_1}^G(S,G_1)
\Psi_{S_1}^{P\cap G_1}(Y,\mu).
\end{multline*}
According to Theorem~\ref{decPsi}, the last line equals
$\Psi_S^P(Y,\mu)$. In the second line, we split the sum over $w$ as
\[
\sum_{\substack{w\in W_0^S\backslash W_0^G\\wL\subset S\\\mu\in\t^*(w\tau^S)}}
\sum_{\substack{w'\in W_0^M\backslash W_0^S\\w'wL\subset M\\\mu\in\t^*(w'w\tau^M)}}
c_{M_1,w'wL}^{P\cap M,U}(w'w\tau,\mu)
\]
and apply Theorem~\ref{decc} with $(G,M,S,L,\tau)$ replaced by
$(S,S_1,M,wL,w\tau)$, which gives
\[
\Phi_{S_1,L}(\exp Y,\tau)
=\sum_{\mu\in\t^*(\tau^G)}\sum_{S\in\CL(S_1)}
\sum_{\substack{w\in W_0^S\backslash W_0^G\\
wL\subset S\\\mu\in\t^*(w\tau^S)}}
c_{S_1,wL}^{P\cap S,U}(w\tau,\mu)\Psi_S^P(Y,\mu).
\]
Changing the order of summation and renaming $S_1$ to~$M$, we obtain the assertion.
\medskip

{\noindent\it Example.} If we specialise to a split maximal torus $T=M=M_0$ in the example from the end of the preceding section, so that $\sigma$ is just a unitary character of~$T(\R)$ with differential~$\mu$, then $c_{T,T}^P(\sigma,w\mu)=m_M(w\check\sigma,P)$ for all $w\in W_T^G=W$, and
\[
\Phi_{T,T}(\exp Y,\sigma)=\sum_{S\in\CL}\sum_{w\in W}
m_M^S(w\check\sigma,P\cap S)\Psi_S^P(Y,w\mu)
\]
for $Y\in\t_P(\R)$ and $\sigma\in \Pi\temp^{P-\mathrm{reg}}(T)$. The version of this result for normalising factors~$r$ was stated as Theorem~7.2 in~\cite{ho-3}, where the contragredience sign on~$\sigma$ was missing.

\bigskip

\noindent\parbox{\textwidth}{\small
Bielefeld University, Germany\\
\texttt{hoffmann@math.uni-bielefeld.de}
}

\end{document}